\documentclass[10pt,twoside]{article}
\usepackage{amssymb, amsmath, amsfonts}
\usepackage{latexsym}
\usepackage{txfonts}
\usepackage{color}
\usepackage{graphicx}
\usepackage{epstopdf}

\begin{document}
\title{ \vspace{0.5in} {\bf\Large A Linking/$S^1$-equivariant Variational
Argument in the Space of Dual Legendrian Curves  and the Proof of
the Weinstein Conjecture on $S^3$ "in the Large"}\\
\vspace{3mm}\normalsize \bf\it To Jalila Ben Othman, in loving
memory } \vspace{3cm}

\author{
{\bf \large Abbas Bahri\vspace{1mm}}\\
{\it\small Department of Mathematics, Rutgers University}\\
{\it\small 110 Frelinghuysen Rd., Piscataway, NJ 08854-8019}\\
{\it\small e-mail: abahri@math.rutgers.edu }} \vspace{1cm}
\maketitle \newpage
\begin{center}
{\bf\small Abstract}\\
\vspace{3mm} \hspace{.05in}\parbox{4.5in} {\small  Let $\alpha$ be a
contact form on $S^3$, let $\xi$ be its Reeb
  vector-field and let $v$ be a non-singular vector-field in $ker
  \alpha$. Let $C_\beta$ be the space of curves $x$ on
  $S^3$ such $\dot x=a\xi+bv, \dot a=0, a \gneq 0$. Let $L^+$,
  respectively $L^-$, be the set of curves in $C_\beta$ such that $b
  \geq 0$, respectively $b \leq 0$. Let, for $x \in C_\beta$, $J(x)=\int_0^1\alpha_x(\dot x)dt$. The framework of the present
  paper has been introduced previously in eg [3].

\hspace{3mm} We establish in this paper that some cycles (an
infinite number of
  them, indexed by odd integers, tending to $\infty$) in the
  $S^1$-equivariant homology of $C_\beta$, {\bf relative} to $L^+ \cup
  L^-$ and to some specially designed "bottom set", see section 4,
   are achieved in the Morse complex of $(J, C_\beta)$ by unions of unstable manifolds of critical
  points (at infinity)which must include periodic orbits of $\xi$;
  ie unions of unstable manifolds of critical points at infinity
  alone cannot achieve these cycles. At the odd indexes $(2k-1)=1+(2k-2)$, $1$ for the linking, $(2k-2)$ for the
  $S^1$-equivariance, we find that the equivariant
  contributions of a critical point at infinity to $L^+$ and to
  $L^-$ are fundamentally asymmetric when compared to those of a
  periodic orbit [5].
The topological argument of existence of a periodic orbit for $\xi$
turns out therefore to be surprisingly close, in spirit, to the
linking/equivariant argument of P.H. Rabinowitz in [12]; eg the
definition of the "bottom sets" of section 4 can be related in part
to the linking part in the argument of [12]. The objects and the
frameworks are strikingly different, but the original proof of [12]
can be recognized in our proof, which uses degree theory, the
Fadell-Rabinowitz index [8] and the fact that
$\pi_{n+1}(S^n)=\mathbb{Z}_2, n\geq 3$. We need of course to prove,
in our framework, that these topological classes cannot be achieved
by critical points at infinity only, periodic orbits of $\xi$
excluded, and this is the fundamental difficulty.

\hspace{3mm} The arguments hold under the basic assumption that no
periodic orbit of index $1$ connects $L^+$ and $L^-$. It therefore
follows from the present work that either a periodic orbit of index
$1$ connects $L^+$ and $L^-$ (as is probably the case for all three
dimensional overtwisted [8] contact forms, see the work of H.Hofer
[10], the periodic orbit found in [10] should be of index $1$ in the
present framework); or (with a flavor of exclusion in either/or) a
linking/equivariant variational argument a la Paul Rabinowitz [12]
can be put to work. Existence of (possibly multiple) periodic orbits
of $\xi$, maybe of high Morse index, follows then.

Therefore, to a certain extent, the present result runs, especially
in the case of three-dimensional overtwisted [8] contact forms,
against the existence of non-trivial algebraic invariants defined by
the periodic orbits of $\xi$ and independent of what  $ker \alpha$
and/or $\alpha$ are. }
\end{center}

{\it \footnotesize Key words}. {\scriptsize  Contact Form,
Legendrian curves, Fadell Rabinowitz index, critical
points at infinity, periodic orbits{\vspace{-1mm}\\
contractible, stable manifold, unstable manifold, dimension,
deformation}}

\medskip

\section{Introduction}

\medskip

Let us consider the Morse relation:

$$\partial c_{2k}^{(\infty)}=c_{2k-1}+h_{2k-1,\infty} (*)$$
 see [1], Lemma 2.14, p126, where $h_{2k-1,\infty}$ is the closure
of a collection of unstable manifolds of critical points at infinity
of Morse index $(2k-1)$ dominated by a collection of periodic orbits
of
 the Reeb vector-field of $\alpha$, $\xi$, of index $2k$, $y_{2k}$s (they can be reduced to a single one, we do not use this here) and where $c_{2k-1}$
 is the closure of a collection of unstable manifolds of periodic orbits
 of $\xi$ of Morse index $(2k-1)$ satisfying the relation $\partial_{per}c_{2k-1}=0$. $\partial_{per}$ is
 the intersection operator related to the periodic orbits of $\xi$. $c_{2k-1}$ is
   assumed here to be a minimal cycle (see [1]), which  means that $c_{2k-1}$ cannot be decomposed
   into smaller cycles for $\partial_{per}$.
 Let $\Gamma_{2s}=\{\text{set of curves made of $s \pm v$-jumps
alternating with $s$ pieces of $\xi$-orbits}\}$

 Let $L^+$ be
 the set of curves in $\cup \Gamma_{2s}$ having all
 their $\pm v$-jumps along $+v$ and let $L^-$ be the set of curves in $\cup \Gamma_{2s}$
  having all their $\pm v$-jumps oriented along $-v$. Let $D_1$ be
  an appropriate
  neighborhood of the critical points (at infinity) of index $1$ of
  $J$, derived by flowing down along the unstable manifolds of these
  critical points small neighborhoods of zero in their stable
  manifolds, see section 4, its figure also for more precisions.

 The first result of this paper states, in a first and rough formulation,
 that the Fadell-Rabinowitz index [8] of the intersection $h_{2k-1, \infty}\cap (J^{-1}([\epsilon, \infty))\smallsetminus D_1)$, is at
 most $(k-2)$\footnotemark\footnotetext{Observe that, unlike in [12]
 and also in [1], we take here for definition of the Fadell-Rabinowitz index of a
 topological
 set $X$ with a free or effective $S^1$-action and classifying map $f$, the power $m$ to which the cohomological Chern class $[x]$ of
  $\mathbb{PC}^\infty$ can be raised and $f^*([x]^m)$ is non zero in the rational cohomology
 of $X$ as in [8]. The Fadell-Rabinowitz index of $\mathbb{PC}^m$ is
 therefore $m$, compare to [12], Lemma 1.13: the Fadell-Rabinowitz
 index of $S^{2m-1}$ for Paul Rabinowitz in [12] is normalized to be $m$, one more
 than we would find with the present definition-which is also the
 definition in [8]-for the quotient $\mathbb{PC}^{m-1}$ of $S^{2m-1}$ by the action of $S^1$.
  We find this definition to be more convenient for our purpose.}.
  The removal of $D_1$ from $J^{-1}([\epsilon, \infty))$ is needed
  in order to warrant that the "bottom set" of $X$, which is $X\cap
  (J^{-1}(\epsilon)\cup \partial D_1)$, is connected in dimension
  $(2k-2)$, since there are no critical points of index $1$ in the
  Morse complex of $X$. We will need to modify this later.

   We then find that the
   proof of the estimate from above on the Fadell-Rabinowitz index of $h_{2k-1, \infty} \cap (J^{-1}([\epsilon, \infty))\smallsetminus D_1)$
 derives from a more general argument: considering a stratified  set $\tilde {X}$, of top dimension $2k$,
 we assume that $\tilde X$ is a manifold in dimensions $2k, (2k-1)$ and that $S^1$ acts effectively
 on $\tilde X$ and freely on its cells of dimension $2k$ and $(2k-1)$ and $(2k-2)$. We also assume
 that we are given an $S^1$-invariant functional $J_\infty$ on $\tilde X$ and a corresponding
 $S^1$-invariant flow such that $\tilde X$ is the closure of the closure of the union of
 unstable manifolds for this flow. We assume that the Palais-Smale condition holds and that $\tilde X$ does not contain any critical point of index $1$.

  Under the above assumptions, we claim that
 that $X=\tilde {X}/S^1$ is of Fadell-Rabinowitz index $(k-2)$ and that there is a classifying
 map for $\tilde {X} \longrightarrow X$ valued into $S^{2k-3} \longrightarrow \mathbb{PC}^{k-2}$.

 Although our argument will contain the proof of the more general claim above, we will
 provide this proof within the framework of Contact Form Geometry [1], [3], [4] and we will discuss
 mainly the case when $X=h_{2k-1, \infty}\cap (J^{-1}([\epsilon, \infty))\smallsetminus D_1)$, that is
 $X=\overline {\cup W_u(y_{2k-1, \infty})}\cap (J^{-1}([\epsilon, \infty))\smallsetminus D_1)$, notations of [1].
 In section $10$ of this paper,
 we will show how the definition of $X$ can be modified in our specific case in order
 to derive the verification of the assumptions above.
 \

 However, this
result does not suffice to impede the above Morse relation since the
same conclusion holds true for the collection of periodic orbits
$c_{2k-1}$ as well, ie $\overline {W_u(c_{2k-1})} \cap
(J^{-1}([\epsilon, \infty))\smallsetminus D_1)$ is also of
Fadell-Rabinowitz index $(k-2)$.

 For $(*)$ above to be impossible, we need a more involved
  estimate on the Fadell-Rabinowitz index of the Morse complexes of dimension $(2k-1)$ {\bf relative} to the values of the
 classifying maps on the topological boundary of these Morse complexes as
 deformation occurs from a collection of periodic orbits $c_{2k-1}$
 to  $h_{2k-1, \infty}\cap (J^{-1}([\epsilon, \infty))\smallsetminus D_1)$, see the
 Morse relation above and see section 11, below, of this paper.

Indeed, the main difference between the case of the periodic orbits
$\overline {W_u(c_{2k-1})}$ and the case of critical points at
infinity $\overline {W_u(h_{2k-1, \infty})}$ stems from the fact
that the periodic orbits "link" the set $L^+$ of curves in
$\cup\Gamma_{2s}$ having only positive $v$-jumps (no $H^1_0$-index
if critical at infinity) with the set $L^-$ of curves in
$\cup\Gamma_{2s}$ having only negative $-v$-jumps (again no
$H^1_0$-index if critical at infinity). This "linking" occurs
because of the first eigenfunction of the linearized operator at a
periodic orbit, see [3], [5].

On the other hand, whereas "linking" of $L^-$ and $L^+$ occurs as a
result of the existence of periodic orbits, at the "bottom level",
in $J^{-1}(\epsilon)$, this linking does not occur and
$J^{-1}(\epsilon) \cap L^+$ and $J^{-1}(\epsilon)\cap L^-$ are
disconnected. They are connected through critical points of index
$1$.

Let $W_1$ be the union of their unstable manifolds (of dimension
$1$). The "linking" induced by the periodic orbits can be recognized
on the classifying maps. Namely, using the map "b" of [5], it is
proven in [5] that  that the pair $(A,B)$, where

$$A=\overline {W_u(c_{2k-1})\smallsetminus
   (L^+\cup L^-)}$$
   and
   $$B= (\overline {W_u(c_{2k-1})\smallsetminus
   (L^+\cup L^-)}) \cap [(\partial L^+ \cup \partial L^-)\cup  {J}_\infty^{-1}(\epsilon)\cup W_1]\cup
   (\overline {\partial_\infty(c_{2k-1}\smallsetminus (L^+ \cup
   L^-)))}$$

   maps through the pair
$$(C_\beta \smallsetminus (L^+\cup L^-), (C_\beta -(L^+ \cup L^-)) \cap
   (\partial (L^+ \cup L^-) \cup {J}^{-1}(\epsilon))\cup W_1\cup
   A_{2k-2})$$, where $A_{2k-2}$ is the set of curves in $C_\beta$
   such that $v$-component $b$ has at least two zeros and at most
   $(2k-2)$ zeros,
   into the pair
$$(\mathbb{P}\mathbb{C}^{\infty}
   \times [-1,1], \mathbb{P}\mathbb{C}^{\infty}
   \times \{-1,0,1\}\cup \mathbb{PC}^{k-2} \times [-1,1])$$
and the composition is onto one of the generators of the homology of
dimension $(2k-1)$ in the target (There are two such generators
since $[-1,1]/\{-1,0,1\}$ has two generators in its homology at
order $1$). $L^+$ and $L^-$ are to be thought in the formulae above
as small attracting (for the decreasing pseudo-gradient)
neighborhoods of these sets.

   On the other hand, each critical point at infinity $h_{2k-1,
   \infty,j}$ in the collection $h_{2k-1, \infty}$ introduces a
   basic asymmetry between $L^+$ and $L^-$, namely $\overline
   {W_u(h_{2k-1, \infty, j})} \cap L^+$ and $\overline
   {W_u(h_{2k-1, \infty, j})} \cap L^-$, one of them, maybe both, is
   of Fadell-Rabinowitz index $(k-2)$ at most, see
   section 7, Lemma 3 below.

  We use this fact and prove that the Morse relation $(*)$ is
   impossible.

   The argument requires some further technical adjustments, which can be completed only under the basic assumption that there are
   no periodic orbit of index $1$.

   Under this assumption,
   we may arrange so that no critical point of index $1$ connects
   $J^{-1}(\epsilon)\cap L^+$ and $J^{-1}(\epsilon)\cap L^-$, see section 4, below.

   The removal of $D_1$ from the sets $X$s above ignores the fact that the periodic orbits link $L^+$ and $L^-$, whereas these two sets are not linked
   in $J^{-1}(\epsilon)$. In order to restore this information, we modify the "bottom set " $J^{-1}(\epsilon) \cup \partial D_1$: we
   "open up" one "side of the bottom set", connecting
   $J^{-1}(\epsilon)\cap L^+$ and $J_0^{-1}(\epsilon)$ (the component
   of $J^{-1}(\epsilon)$
   close to "small" back and forth runs along $v$) and we create in this way
   a new "bottom set" $D_1^+$. $D_1^+ \cup (J^{-1}(\epsilon)\cap L^-)$
   may be viewed, after a re-parametrization of flow-lines and after
   a related definition
   of a new functional $\tilde J$, see J.Milnor [11], Theorem 4.1, pp37-38, as $\tilde
   J^{-1}(\epsilon)$. We now have a disconnected "bottom set" $\tilde J^{-1}(\epsilon)$, where $L^+$ and $L^-$ are not connected
   anymore, but $L^+$ and $J_0^{-1}(\epsilon)$ are connected.

   Let $W_1^-$ be the part of $W_1$ related to
   $L^-$, ie connecting the various components of $J^{-1}(\epsilon)\cap L^-$ exclusively.

   Replacing $J$ by $\tilde J$ in the pairs above and $W_1$ by
   $W_1^-$, we find now a classifying map valued into $(\mathbb{P}\mathbb{C}^{\infty}
   \times [-1,1], \mathbb{P}\mathbb{C}^{\infty}
   \times (\{-1\} \cup[0,1])\cup \mathbb{PC}^{k-2} \times [-1,1])$.
   This pair has the advantage when compared to the former one that
   it has only one generator in dimension $(2k-1)$.

  We can now use the asymmetry of $L^+$ and $L^-$ for $h_{2k-1,
   \infty}$ as described above and we prove, after careful
   modifications that are embedded in an isotopy of decreasing
   pseudo-gradients for the functional, that the Fadell-Rabinowitz
   index of $\overline {W_u(h_{2k-1, \infty})} \cap \tilde
   J^{-1}([\epsilon, \infty))$ is $(k-2)$ {\bf relative} to the
   value of the classifying map on the "bottom set" $B_0=D_1^+ \cup (J^{-1}(\epsilon)\cap
   L^-)\cup W_1^-$, which is constrained to take values into $\mathbb{P}\mathbb{C}^{\infty}
   \times \{-1\} \cup[0,1])\cup \mathbb{PC}^{k-2} \times [-1,1]$.

   The same conclusion cannot hold for $\overline {W_u(c_{2k-1})} \cap \tilde
   J^{-1}([\epsilon, \infty))$ and the contradiction argument
   follows.

   As stated above, a basic assumption is used in this argument:
   namely, it is assumed that the sets $J^{-1}(\epsilon)\cap L^+$ and $J^{-1}(\epsilon)\cap L^-$
   are not connected by a periodic orbit of index $1$.

   We conjecture that, in the framework of over-twisted contact
   forms [7], the periodic orbit found by H.Hofer [10] is of index $1$
   (when viewed in our framework).

   In some regards, our present paper indicates that, for the
   existence of periodic orbits, either an equivariant/linking
   argument "a la Paul Rabinowitz" [12] works, yielding a sequence of
   periodic orbits of odd Morse index $(2k-1)$ for $k$ large; or
   this argument does not work and a periodic orbit of index $1$ is
   found, as in H.Hofer [10] (maybe and probably).

   This is not established rigorously, but strongly indicated by the
   proof. This is emphasized in the last section of this paper.

     Theorem 1.3,(i) of [1], the proof of which was not complete, see [2],
  follows from the claims above:\\
  \medskip

\noindent{\bf Theorem 1}\, {\it Assume that $\alpha$ is a contact
form on $S^3$ and that the Reeb vector-field of $\alpha$ has no
periodic orbit
  of Morse index $1$. Then, (*) is impossible for $k$ large enough and $J$ has a sequence of critical values corresponding to periodic orbits
  of index $(2k-1)$.}

  \medskip
  Let us recall that the existence of one periodic orbit for the contact forms of the tight contact structure
 of $S^3$ is a theorem by P.H.Rabinowitz [12], established without dimension restriction, whereas
 the existence of one periodic orbit for the
 contact forms of all over-twisted [7] contact structures on a closed contact
 three dimensional manifold is a theorem by H.Hofer [10].

  Theorem $1$ above gives a new proof for the Weinstein conjecture
 on $S^3$. This new proof combines the case of the tight contact structure on $S^3$ and the case of all the other over-twisted ones [7] and, therefore,
 could lead to a better understanding of the existence process for periodic orbits of $\xi$. This new proof could also possibly lead to multiplicity results,
  on all three dimensional closed
  contact manifolds with finite fundamental group.

  The present paper and the
 corresponding
 topological argument for existence show also how to overcome the non-compactness of the variational
 problem associated to the periodic orbits problem for
 the Reeb vector-field $\xi$ of a given contact form $\alpha$ on a three dimensional
 closed contact manifold with finite fundamental group.\\

\noindent  {\bf 2.  The Fadell-Rabinowitz index of\\ $X=\overline {\cup W_u(y_{2k-1, \infty})}\cap (J^{-1}([\epsilon, \infty))\smallsetminus D_1)$}\\

\noindent By assumption, the set $X$ can be written as a union of
closures of unstable manifolds of critical points at infinity of
index $(2k-1)$

$$X= \overline {\cup W_u(y_{2k-1, \infty})}\cap (J^{-1}([\epsilon, \infty))\smallsetminus
D_1)$$

As stated above, $D_1$ is derived after flowing down along the
unstable manifolds of dimension $1$ of the critical points (at
infinity) of index $1$ of $J$ small neighborhoods of zero
(transverse to the flow) in their stable manifolds, see section 4 in
order to recognize this construction with the help of a drawing.

Let us assume that $X$ is a manifold in dimensions $(2k-1)$ and
$(2k-2)$, see section 10 for the verification of these assumptions.
It follows from these assumptions that each $y_{2k-1, \infty}$  is
simple and that there cannot be more than one flow-line from each to
$y_{2k-1, \infty}$ to each $y_{2k-2, \infty}$. This observation
helps the understanding. We then claim:\\

\noindent {\bf Lemma 1}  {\it The Fadell-Rabinowitz index of $X$ is
$(k-2)$ and there is a classifying map for the $S^1$-action on $X$
valued into $S^{2k-3}/\mathbb{PC}^{k-2}$.}

\medskip
\noindent {\it Proof of Lemma 1.}  Let us consider the topological
boundary of each $\overline { W_u(y_{2k-1, \infty})}\cap
(J^{-1}([\epsilon, \infty))\smallsetminus D_1)$, which we denote
$Z_{2k-2, \infty}$. It is a chain of dimension $(2k-2)$. Let

$$f: Z_{2k-2, \infty} \rightarrow \mathbb{PC}^\infty$$
be any classifying map for the $S^1$-action on $\tilde {Z}_{2k-2,
\infty}\rightarrow Z_{2k-2, \infty}$, where $\tilde {Z}_{2k-2,
\infty}$ is the set of $S^1$-invariant curves over $Z_{2k-2,
\infty}$, $\tilde {Z}_{2k-2, \infty}=S^1 *Z_{2k-2, \infty}$.

We may assume $f$ to be $C^\infty$, so that, by general position,
its image may be assumed, after deformation, to be valued into
$\mathbb{PC}^{k-1}$:

$$f: Z_{2k-2, \infty} \rightarrow \mathbb{PC}^{k-1}$$

Using then degree theory, we may assume that $deg f=0$, since
$Z_{2k-2, \infty}$ is a boundary. Observe that $Z_{2k-2, \infty}$ is
connected, being the image through the time $1$-map of the
decreasing pseudo-gradient acting on an unstable sphere $S^{2k-2}$
for $y_{2k-1, \infty}$.

In the special framework of [1], [3], [4] and [5], with $C_\beta=\{x
\in H^1(S^1,M); \beta(\dot x)=d\alpha(\dot x,v)=0, \alpha(\dot
x)=\text {C;C not prescribed and positive}\}$, $v$ non singular in
$ker\alpha$ and with $\cup \Gamma_{2s}$, we may also assume that $f$
is given, on the set of curves $x$ such that the $v$-component $b$
of their tangent vector $\dot x$ has at least two zeros and at most
$(2k-2)$ zeros and is equal to the map "$b$" of [5], which, after
deformation, is then valued in $\mathbb{PC}^{k-2}$; this in the
specific case as in [3], [4], [5] etc. In other cases, $f$ might be
given on some other set that maps into $\mathbb{PC}^{k-2}$. For
simplicity and in order to make our arguments more transparent, we
assume in the sequel that we are in the specific framework of [1],
[3], [4], [5]. The generalization is clear.

Let us pick up a point $x_0$, which is a regular value for $f$ in
$\mathbb{PC}^{k-1}$, not in $\mathbb{PC}^{k-2}$ and let us consider
$f^{-1}\{x_0\}$. If there are no points in $f^{-1}\{x_0\}$, our
argument is complete, see below. Otherwise, we assume for sake of
simplicity that $f^{-1}\{x_0\}=\{z_1, z_2\}$, that is it is made of
exactly two points where $f$ has Jacobians with opposite signs.

We then consider a generic path from $z_1$ to $z_2$. We can choose
$x_0$ so that $z_1$ and $z_2$ are not in the stable manifold of any
critical point (at infinity) of $X$. Using then the decreasing
pseudo-gradient that defines $X$, we can deform this path into a
path in $W_{2k} \cup L^+ \cup L^-\cup J_0^{-1}(\epsilon) \cup D_1$.
$W_{2k}$ is the set of curves such that the $v$-component $b$ of
their tangent vector has $2k$-sign changes, not less-we may assume
that if $z_i$ has $2k$ zeros, then these $2k$ zeros survive all
along the decreasing flow-line, until the "bottom set" is reached,
this is not essential in the argument, it is rather a side remark-,
$L^+$ is the set of curves where $b$ is positive, $L^-$ the set with
$b$ negative, $J_0^{-1} (\epsilon)$ is the component of
$J^{-1}(\epsilon)$ made of "small" curves in $J^{-1}([0,
\epsilon])$, close to back and forth or forth and back runs (one or
several) along $v$ and $D_1$ is a small neighborhood of the unstable
manifolds of the critical points (at infinity) of $J$ of index $1$
deleted from a small neighborhood of its trace in $L^+ \cup L^-$.
For this reason, all the curves of $\partial D_1 \smallsetminus
(L^+\cup L^-)$  are such that their $v$-component $b$ has at least
two zeros (see Lemma 7 for more precisions in this specific case,
the argument extends to the general case, when the classifying map
is not $"b"$ anymore). Since the Morse index of these critical
points (at infinity) is $1$, we find that the union of the unstable
manifolds of these critical points at (infinity)is a compact set and
we find that $b$ on this neighborhood can be deformed to a function
$\tilde b$ having a finite, a priori bounded number of zeros, given
by the projection of $b$ onto the unstable directions, so that
$\partial D_1\smallsetminus (L^+\cup L^-)$ can be mapped through a
modification of the map $b$ in to $\mathbb{PC}^r$, for a fixed $r$
independent from $k$.

Let $B_1=J_0^{-1}(\epsilon)\cup [J^{-1}(\epsilon)\cap (L^+ \cup
L^-)] \cup \partial D_1$

We use this path and standard methods, see M.Hirsch [9], pp126-127
and we modify the map $f$ near the path and make it valued into
$\mathbb{PC}^{k-1}\smallsetminus \{x_0\} \cong \mathbb{PC}^{k-2}$.
Let us outline in details the argument:

Let $p$ be the path as above. After deformation,
 we may assume that this path takes the following form: $p$ starts at $z_1$ with a decreasing flow-line in the corresponding $W_u(y_{2k-2, \infty})$.
 This flow-line will, using general position, reach the "bottom set" $B_1$; same for $z_2$,
 and this happens whereas the flow-lines do not leave their respective $W_u(y_{2k-2, \infty})$s.
 There are no critical points of index $1$ above $B_1$ by construction and therefore,
 we may assume that the remainder of the path $p$ is in a subset $Z$ which is a manifold in dimension $(2k-2)$ and in dimension $(2k-3)$,
 so that $p$ does not cross any singularity in $W_u(y_{2k-1, \infty})$. The cancellation procedure of section $1$ and of [9], pp126-127, may be applied.
By general position, we can assume, for a given copy of
$\mathbb{PC}^{k-2}$, that $f(p)\cap \mathbb{PC}^{k-2}=\varnothing$.
Thus, we may assume that $f(p)$ and in fact $f(D^{2k-2})$ is
contained in a disk $D_2^{2k-2}$ around $x_0$.

 As a map,  from $D^{2k-2}$ into $D_2^{2k-2}$, $f$ is then, using a
degree argument, homotopic relative to its boundary value (from
$\partial D^{2k-2}$ into $\partial D_2^{2k-2}$) to a map valued into
$\partial D_2^{2k-2}$. Using an equivariant family of small sections
to the $S^1$-action in $S^{2k-1}$ and using the lift $\hat f$ of
$f$, we can lift this homotopy into a homotopy of $S^1$-equivariant
maps above. Since $\hat f(\tau * x)=e^{ip\tau}f(x)$, the same
relation will hold for all lifts along the homotopy and, at the end
of this homotopy, the classifying map for $\overline {\partial
(W_u(y_{2k-1, \infty}) \cap J^{-1}([\epsilon, \infty))}$ will be
valued into $S^{2k-1}\smallsetminus S^1 *\{x_0\}$, thus it will be
valued into $S^{2k-3}$, $\mathbb{PC}^{k-2}$ as claimed, with the map
unchanged on the set of curves where $b$ has at least one
sign-change and at most $(2k-2)$ zeros, as claimed.

We find then a new map $\tilde f$, {\bf equal} to $f$ on the set
where $b$ has at least one zero (with a sign change) and at most
$(2k-2)$ zeros.

We extend now this map, or rather some power of this
 map to all of $\overline { W_u(y_{2k-1, \infty})}\cap \tilde {J}^{-1}([\epsilon, \infty))$.
 $\tilde f$ can be assumed to be defined in fact on all of $\overline {W_u(S^{2k-2})} \cap  (J^{-1}([\epsilon, \infty)\smallsetminus D_1))$,
 since this set retracts by deformation on $\partial\overline { W_u(y_{2k-1, \infty})}\cap (J^{-1}([\epsilon, \infty)\smallsetminus D_1))$. Restricting,
 it follows that $\tilde f$ is defined from $S^{2k-2}$ into $\mathbb{PC}^{k-2}$. Lifting,
 we find an equivariant map $\hat f: S^{2k-2}\rightarrow S^{2k-3}$. $\hat f$ is equivariant
 in that $\hat f(e^{i\tau}*x)=e^{pi\tau}\hat f(x)$, for a given integer $p$. This is with an
 appropriate modification of the map $b$, where the various component of $b$ on the various
 functions $sin(2j\pi t)$ and $cos (2j \pi t)$ are raised at the appropriate powers so that
 the modified map, with the introduction of this powers, satisfies the equivariant law as
 written above, see [5] for the transformation of $b$ into its $L^2$-projection on the appropriate Fourier modes.

Let us restrict the map $\hat f$ to $S^{2k-2} \times \{1\}$, we find
a map $g: S^{2k-2} \rightarrow S^{2k-3}$. We know that the homotopy
group of order $(2k-2)$ of $S^{2k-3}$ is $Z_2$ for $k \geq 3$.
Therefore, if we knew that $g$ was a double, we could extend it to
$D^{2k-1}$, thereby extending $\tilde f$, valued into
$\mathbb{PC}^{k-2}$, to all of $\overline { W_u(y_{2k-1,
\infty})}\cap \tilde {J}^{-1}([\epsilon, \infty))$.

In order to be sure that $g$ is a double, we need to be able to
compose it with a map of degree $2$, or a map of even degree from
 $S^{2k-3}$ into itself. There are such maps and, thinking in terms of
 the covering map $h:S^1*\overline { W_u(y_{2k-1, \infty})}\cap (J^{-1}([\epsilon, \infty)\smallsetminus D_1)) \rightarrow S^{2k-3}$
 over $\tilde f$, we can assume that, in order to define $h$, we have composed its original value
 with a map from $S^{2k-3}$ into itself and we have raised each (complex) component
 to the power $2$ and re-normalized thereafter so that the norm stays $1$. The resulting
 map is equivariant: it does satisfy the law $h(e^{i\tau}*x)=e^{pi\tau}h(x)$ with a
 suitable $h$, for which there is a suitable $p$. After this composition, the map $g$ that
 we find is equal to the previous value for $g$ composed with a map of even degree from $S^{2k-3}$ in
 itself and it follows that the new map $g$ is a double and the extension can be completed.

In this way, we find that the map "$b$", defined on the set of
curves having
 a least one sign-change and at most $(2k-2)$ zeros, appropriately modified
 by reducing it to its orthogonal $L^2$-projection on the basis of functions
 $sin (2j\pi t), cos (2j\pi t), 1\leq j \leq (k-1)$, also appropriately
  modified by raising these components to the appropriate powers and by
  taking only "part" of this map on the $U_1$ as above, that is changing
  $b$ on $U_1$  into its projection on the corresponding negative eigenfunction(s)
  (and thereby finding a function valued in a finite dimensional fixed
  $\mathbb{C}^{r+1}$), we find that this modified map "$b$" extends to
  all of $h_{2k-1, \infty}$ into a map which is equivariant with the use of an $e^{ip\tau}$ factor of covariance in lieu of
  $e^{i\tau}$. The claim follows. $\blacksquare$\\

\noindent {\bf 3. $h_{2k-1, \infty}$ and $c_{2k-1}$, splitting of
the
argument above and introduction of a Basic Assumption}\\

\noindent The above argument is insensitive to the fact that the
$y_{2k-1, \infty}$s are periodic orbits or critical points at
infinity. This is essentially due to the fact that the "bottom set"
$B_1$ is "above" any critical point of index $1$, so that $L^+$ and
$L^-$ can be connected through this "bottom set". We need, in order
to distinguish between the case of the periodic orbits and the case
of the critical points at infinity, to keep $L^+$ and $L^-$
separated in the "bottom set".

We are therefore led to introduce the following basic assumption in
our work:

  (A) $L^+$ and $L^-$ are not connected by
  a periodic orbit of index $1$.

   We also assume that each of $L^+\cap J^{-1}(\epsilon)$ and $L^-\cap J^{-1}(\epsilon)$
   is connected to the
  "small" (these are the curves of $C_\beta$ close to one or several back and forth or forth and back runs along $v$, they are
  contractible in a given, small neighborhood of eg their base point) curves of $J^{-1}(\epsilon)$ by a critical
  point of index $1$, respectively $x^{1, \infty}_+$ and $x^{1,
  \infty}_-$. After re-parametrization of the flow-lines of a pseudo-gradient for $J$ which modifies this functional, but leaves
  the flow-lines of the pseudo-gradient unchanged see J.Milnor [11], Theorem 4.1, pp37-38, and after tangencies between critical
  points of index $1$, we may assume
  that these are the only critical points of index $1$ connecting the
  "small"
  contractible curves (as above) of $C_\beta$ to $L^+$ and to $L^-$.  Using
  this re-parametrization procedure [11] and again tangencies, we may
  also assume then that $L^+$ and $L^-$ are not connected by
  critical points at infinity of index $1$: The unstable manifold of such a critical point at infinity $\bar x^\infty$,
  on the side going to $L^+$ or on the side going to $L^-$, is made of curves having changes in the orientations of their $\pm v$-jumps.
    This change of sign allows, without disturbing the flow-lines in $L^+$ and in $L^-$, to complete a tangency (maybe after
  re-parametrizing the flow-lines and changing the functional as in J. Milnor [11]) with $x^{1, \infty}_+$ or with
  $x^{1, \infty}_- $ and remove the direct connection  between  $L^+$ and  $L^-$.
$ L^+$ and $ L^-$-we might need to change $J$ into $\tilde J$- are
then
  not anymore
  directly connected by critical points of index $1$. They are
  connected through the "small" contractible curves of $C_\beta$.

  All the re-parametrizations and tangencies completed above do not
  perturb the flow-lines in $L^-$ and in $L^+$.

  The most general form of our basic assumption is that we do not
  have a periodic orbit of index $1$ connecting curves of $L^+$ with
  curves of $L^-$ whereas there would be at the same time critical
  points at infinity of index $1$ connecting $L^+\cap J^{-1}(\epsilon)$ and the "small"
  contractible (as above, in a given small neighborhood of eg their base point) curves of $C_\beta$ and connecting $L^-\cap J^{-1}(\epsilon)$ and the
  "small" contractible curves of $C_\beta$. If this assumption does not hold, we would find a
  "circle" of critical points of index $1$ between $L^+$, $L^-$ and
  the "small" contractible curves and our arguments then collapse.
  As long as some separation occurs along this circle, it appears
  that the above arguments goes through.\\

\noindent {\bf 4. Bottom Sets}\\

   \noindent Our "bottom set" $B_1$ above, which is $J^{-1}(\epsilon) \cup
  \partial D_1$, is connected. This does not allow to recognize the
  contribution of the periodic orbits, as described above. We
  therefore define below another "bottom set" $B_0$. In its manifold
  part (outside the unstable manifold of the critical point $x^{1, \infty}_-)$, it disconnects $L^+ \cup J_0^{-1}(\epsilon)$ and
  $L^-$. This of course destroys an essential feature of
our argument above about the Fadell-Rabinowitz index of $X$, namely
that the flow-lines out of $z_1$ and $z_2$ can be connected in the
"bottom set". We cannot assert this anymore with $B_0$. We will see
how to overcome this difficulty.

  We need in fact to define for the purpose of our argument below two distinct
  "bottom sets", $D_1^+$ and $D_1^-$ which are built from the same
  principle, but are different and not symmetric in their
  definition.

  The basic pieces for the definition of $D_1^+$ are
  $J^{-1}(\epsilon) \cap L^+$ and $J_0^{-1}(\epsilon)$, where
  $J_0^{-1}(\epsilon)$ is the component of $J^{-1}(\epsilon)$ made
  of "small" contractible curves of $C_\beta$ (near back and forth or forth and back runs along $v$). These various pieces
  are glued with boundaries of neighborhoods of unstable manifolds
  of the various critical points at infinity of index $1$
  connecting the various components of $J^{-1}(\epsilon) \cap L^+$
  and connecting a component of this latter set to
  $J_0^{-1}(\epsilon)$. Flowing down the boundary (transverse to the flow) of a small neighborhood of $0$ in the the stable manifold of each of
  this critical point of index $1$ on each side of its unstable
  manifold and glueing with the corresponding bottom components of
  $J^{-1}(\epsilon)$ (this requires deletion of a neighborhood of the trace of this unstable manifold on the
  bottom component and glueing them, see the two figures below), we find for $D_1^+$ a manifold which acts
  exactly as a level surface for $J$, ie the flow of a decreasing
  pseudo-gradient is transverse to $D_1^+$.

  For $D_1^-$, we complete the same construction with $J^{-1}(\epsilon) \cap
  L^-$ {\bf only}; that is we do not add $J_0^{-1}(\epsilon)$ and do
  not connect it through the unstable manifold of $x^{1, \infty}_-$
  to $J_0^{-1}(\epsilon)$.
  \medskip
\begin{figure}[!htb]
 \centering
 \includegraphics[scale=.4]{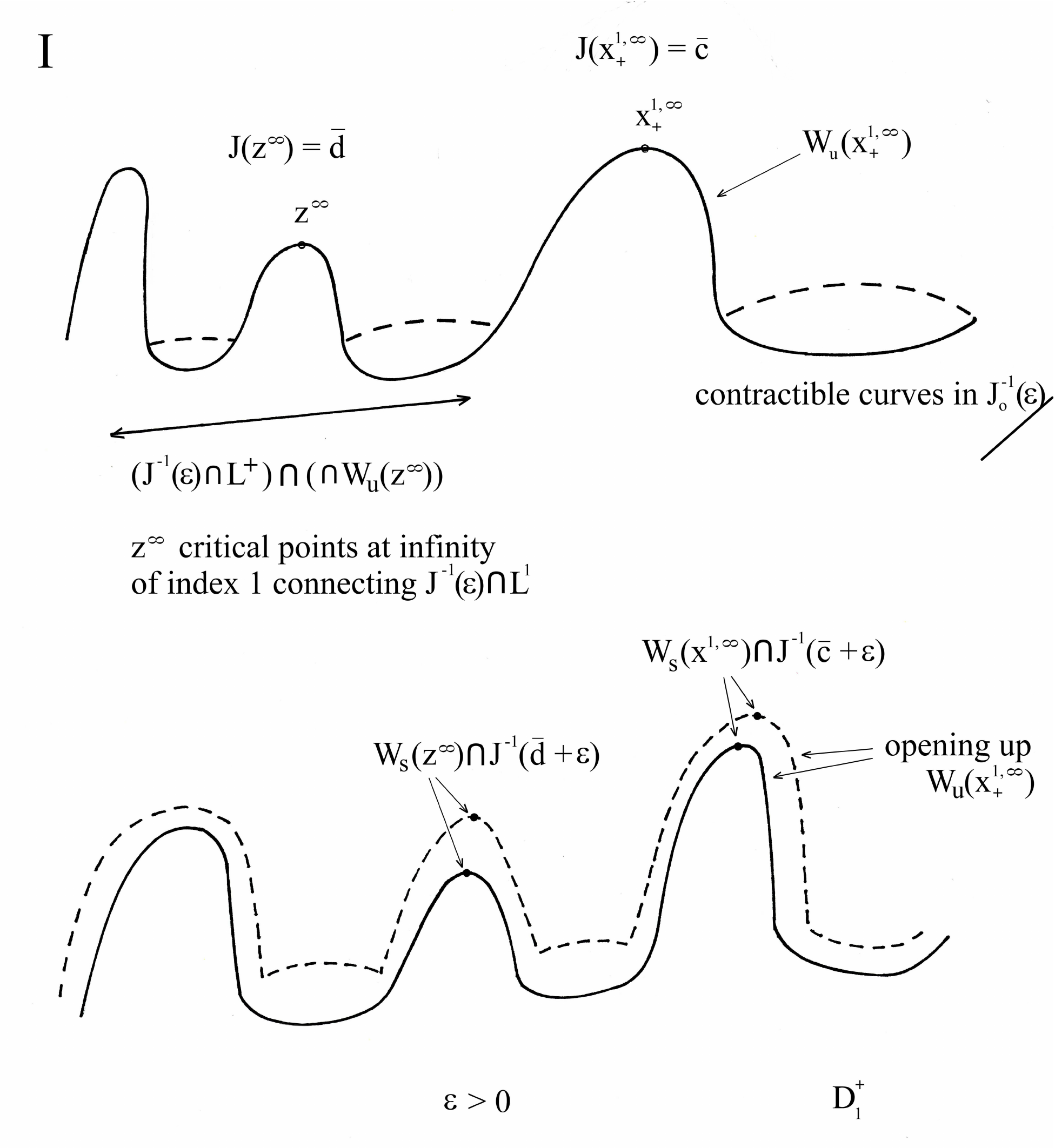}
 \end{figure}

\begin{figure}[!htb]
 \centering
 \includegraphics[scale=.3]{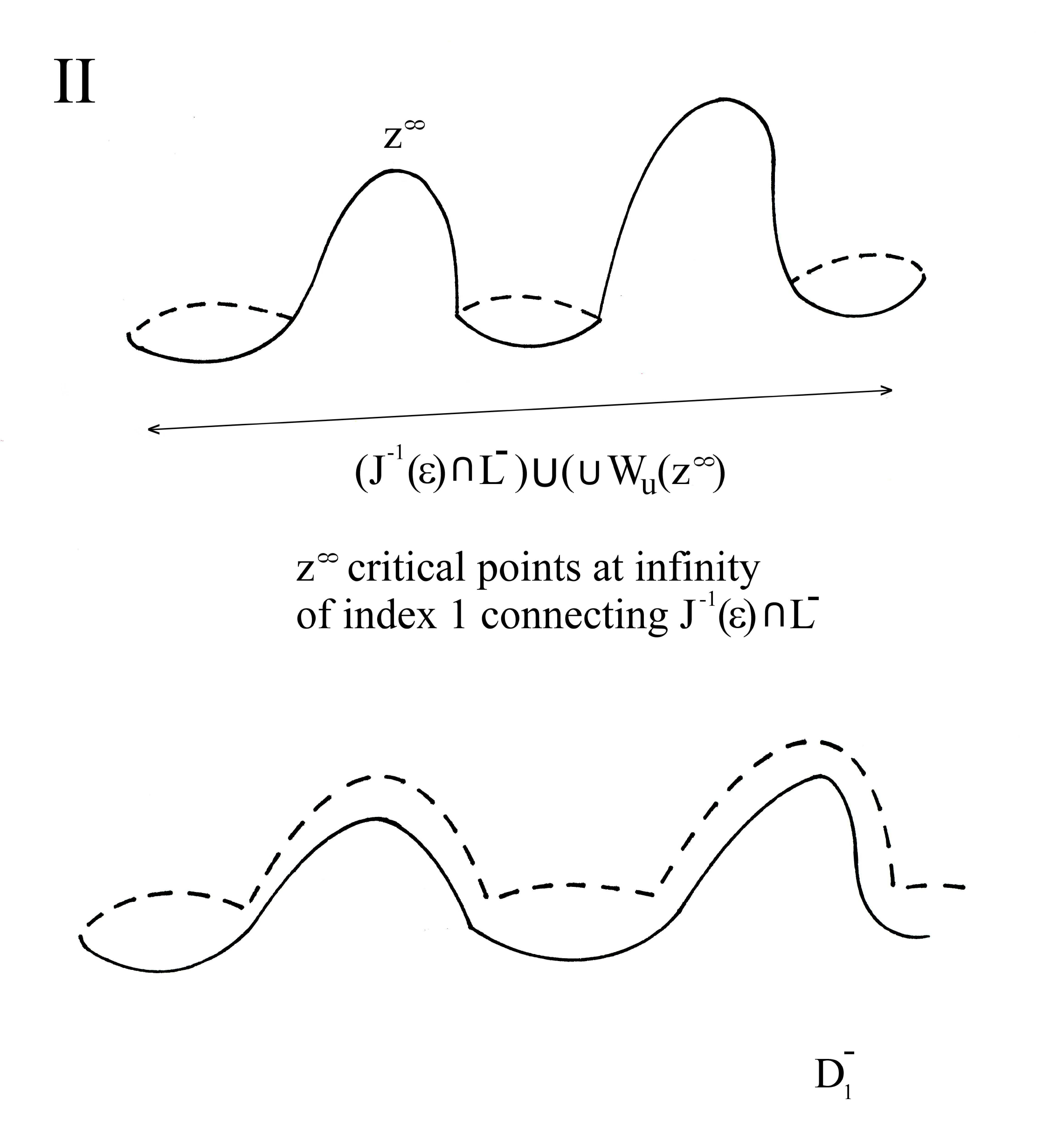}
 \end{figure}


\newpage

 There is a fundamental asymmetry between the definition of $D_1^+$
  and the definition of $D_1^-$.

  For the purpose of our argument, we will denote $U_1$ the part of
  $D_1^+$ which has been built using the stable manifold of the
  critical points of index $1$ connecting $L^+\cap J^{-1}(\epsilon)$ and
  $J_0^{-1}(\epsilon)$ on one hand and connecting the various components of
  $J^{-1}(\epsilon)\cap L^+$ between themselves on the other hand.
  We will denote $B_0$ the union $D_1^+ \cup D_1^-\cup W_u(x^{1,
  \infty}_-)$. The manifold part of $B_0$ is $D_1^+ \cup D_1^-$,
  which disconnects $L^+ \cup J_0^{-1}(\epsilon)$ and $L^-$. This is
  what we have sought.

  As noted above, we may assume that we did re-parametrize the flow-lines of a/the
  decreasing pseudo-gradient just as in J.Milnor [11], Theorem 4.1, pp37-38 and that we thus have derived a new
  functional $\tilde J$ that has the same critical points (at
  infinity) than $J$, with the same stable and unstable manifolds
  for each of these critical points (at infinity) and for which
  $D_1^+ \cup D_1^-$ is $\tilde J^{-1}(\epsilon)$.\\

\noindent {\bf 5. Splitting of the critical points at infinity of
  $h_{2k-1, \infty}$ into two groups}\\

 \noindent  We split the critical points at infinity composing $h_{2k-1, \infty}$ into two groups. In the first group, the $y_{2k-1,j}^
  \infty$ are such that one of their large $\pm v$-jumps is along $+v$, whereas, in the second group, all the large
  $\pm v$-jumps of the critical points at infinity are along $-v$. Completing tangencies, we may
  assume that the second group is reduced to a single
  $z^{\infty,-}_{2k-1}$. We will have to recall, in section 9, that we
  reached this single $z^{\infty, -}_{2k-1}$ out of several such
  critical points at infinity, all of which have their large $\pm
  v$-jumps along $-v$.
  \medskip

\begin{figure}[!htb]
 \centering
 \includegraphics[scale=.35]{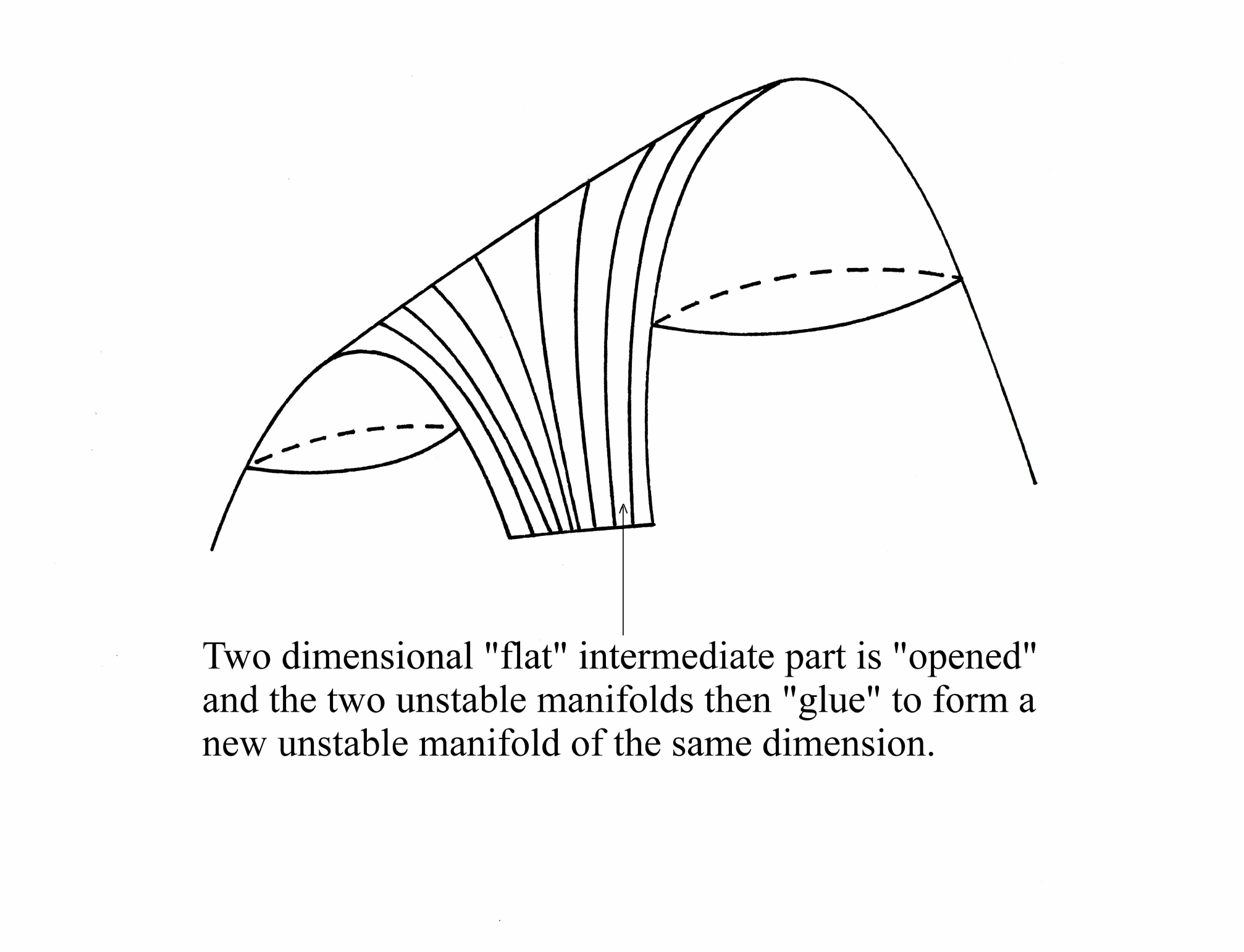}
 \end{figure}

 \noindent {\bf 6. Requirements for the application of the arguments of
  section after the definition of a new "bottom set"
  $B_0$}\\

  \noindent  In order to apply the arguments of section, we now need to
  know that
the traces of $W_u(z_{2k-1,-}^\infty)$ and the trace of each
$W_u(y^{\infty}_{2k-1, j})$ on the components $D^+_1$ and $D^-_1$ of
the bottom set $B_0$ are connected, see section 7 and section 8
below. We also need to know that the trace of
$W_u(z_{2k-1,-}^\infty)$  on $D_1^-$ is connected on the other hand.
These results are established in the next section, after appropriate
modifications of the pseudo-gradient.\\

 \noindent {\bf 7. Preliminary Technical Results}\\

\noindent We start with:\\

\noindent {\it  The classifying map on  $h_{2k-1, \infty}\cap
  J_0^{-1}(\epsilon)$} :\\

\noindent Let $J_0^{-1}(\epsilon)$ be the component of
$J^{-1}(\epsilon)$ corresponding to curves close to back and forth
or forth and back runs along $v$, which we have also have been
referring to as the component of $J^{-1}(\epsilon)$ made of "small"
contractible curves.

We first modify $W_u(y^{\infty}_{2k-1, j })$ with the addition of
"bridges" in order to render $W_u(h_{2k-1, \infty}) \cap
J_0^{-1}(\epsilon)$ connected. This is completed with the
introduction of additional critical points
  $z_{2k-1,j}^\infty$s, of critical value eg $2\epsilon$, which have their boundaries made of
  flow-lines all abutting to "small" contractible curves. Each $z_{2k-1,j}^\infty$
  has in its boundary two companion critical points at infinity of index $(2k-2)$, $ z^i_{2k-2}$, $i=1,2$, which,
  together with $z_{2k-1, j}^\infty$ help build the "bridge. The critical values of these latter points are eg
  $3\epsilon/2$. The functional $J$ is again slightly perturbed, we keep the same notation $J$ or
  $\tilde J$. With these "bridges",  $W_u(h_{2k-1, \infty}) \cap
  J_0^{-1}(\epsilon)$ is now
connected in dimension $(2k-2)$. We then claim that:
\medskip

\noindent {\bf Lemma 2} \,\, {\it The Fadell-Rabinowitz index of
$h_{2k-1, \infty}\cap
  J_0^{-1}(\epsilon)$ is $(k-2)$ at most. After possible addition of
  "bridges", the classifying map for the $S^1$-action on
  $\overline{W_u(h_{2k-1, \infty})}\cap J_0^{-1}(\epsilon)$ may be
  assumed to be valued in $S^{2k-3}/\mathbb{PC}^{k-2}$}
  \medskip

\noindent{\it Proof of lemma 2. } Here $J_0^{-1}(\epsilon)$
designates the level surface $\epsilon$
  of the functional $J$, in the connected component corresponding to
  contractible curves.

The proof of the Lemma starts with the relation:

  $$\partial c_{2k}^{(\infty)}=c_{2k-1}+h_{2k-1,\infty}$$
where $c_{2k-1}$ and $h_{2k-1, \infty}$, as well as
  $c_{2k}^{(\infty)}$ designate the collection of unstable manifolds
  (with closures) of the various critical pints (at infinity)
  involved in the definition of each piece. It then follows that:

  $$\partial c_{2k}^{(\infty)}\cap J_0^{-1}(\epsilon)=\partial (c_{2k}^{(\infty)}\cap
  J_0^{-1}(\epsilon))=c_{2k-1}\cap
  J_0^{-1}(\epsilon)+h_{2k-1,\infty}\cap
  J_0^{-1}(\epsilon)$$
  Since $c_{2k-1}\cap
  J_0^{-1}(\epsilon)$ has a classifying map valued into $S^{2k-3}$,
  the same can be inferred of $h_{2k-1,\infty}\cap
  J_0^{-1}(\epsilon)$ if this set is connected. If not, we have
  resolved this connectedness issue with the addition of a finite
  family of paths, with tubular neighborhoods (following appropriate
  constructions).

  In all, after some required modifications, we may assume that the
  classifying map for every
  trace on $J_0^{-1}(\epsilon)$ of the closure of a collection of unstable manifolds of
  dimension $(2k-1)$, which we assume to be a manifold in dimensions
  $(2k-1)$ and $(2k-2)$, cobordant to $c_{2k-1}$ is valued into
  $S^{2k-3}/\mathbb{PC}^{k-2}$.

  We may add to $c_{2k-1}\cap
  J_0^{-1}(\epsilon)$ and to $h_{2k-1,\infty}\cap
  J_0^{-1}(\epsilon)$ the unstable manifolds of the critical points
  (at infinity) of index $1$ and also  $J^{-1}(\epsilon) \cap (L^+
  \cup L^-)$. Since this latter set is of low Fadell-Rabinowitz
  index, we can assert that the Fadell-Rabinowitz index of the union
  is at most $(k-2)$.

   Recalling our construction above now, when we were defining the "bottom sets", we take the ""side of $L^+$ and "open-up" the unstable
   manifolds
  of dimension $1$ connecting $J^{-1}_0(\epsilon)$ to
  $J^{-1}(\epsilon)\cap L^+$ and connecting the various components of $J^{-1}(\epsilon)\cap L^+$ between themselves, in order to create a "level surface"
  $D_1^+$ transverse to the flow.

  The "opening-up" is completed with the use of the Morse Lemma at
  {\color{red} $x^{1, \infty}_+$ and the various other critical points } at infinity of index $1$ related to $J^{-1}(\epsilon \cap L^+$.
  The "top" of $D_1^+$ at $x_+^{1(\infty)}$ is made
  of the trace of $W_s(x_+^{1(\infty)})$, the stable manifold of
  $x_+^{1(\infty)}$, on a level surface just above
  $x_+^{1(\infty)}$. A neighborhood of this "top" is "flown down" on
  both "sides" of $x_+^{1(\infty)}$ and connects
  $J_0^{-1}(\epsilon)$ and $J^{-1}(\epsilon) \cap L^+$. $D_1^+$ is
  the union of the three pieces $J_0^{-1}(\epsilon)$, $J^{-1}(\epsilon) \cap
  L^+$ and the piece related to these unstable manifolds of
  dimension $1$.

  It is then clear that the Fadell-Rabinowitz index of $c_{2k-1}
  \cap D_1^+$ and of $h_{2k-1, \infty}\cap D_1^+$, as well as that of
  their union, is at most $(k-2)$ since these sets can be
  equivariantly mapped into  $(c_{2k-1}\cap
  J_0^{-1}(\epsilon))\cup W_u(x_+^{1, (\infty)}) \cap (c_{2k-1}\cap
  J^{-1}(\epsilon)\cap L^+)$ and into $(h_{2k-1, \infty}\cap
  J_0^{-1}(\epsilon))\cup W_u(x_+^{1, (\infty)}) \cap (h_{2k-1, \infty}\cap
  J^{-1}(\epsilon)\cap L^+)$ as well as into their union. $\blacksquare$\\

  \medskip
 \noindent{\bf Lemma 3}\,  {\it Let $z^\infty_{2k-1}$ be a critical point at
  infinity of index $(2k-1)$. Let $\partial$ be the
  intersection operator. Then,
 $\partial z^\infty_{2k-1} \cap L^+$ or $\partial
  z^\infty \cap L^-$ is empty for a suitable globally defined,
  admissible (ie leaving $L^+$ and $L^-$ invariant) decreasing
  pseudo-gradient. In fact, the classifying map for the $S^1$-action
  on either $\overline {W_u(z^\infty_{2k-1})} \cap L^+$ or on  $\overline {W_u(z^\infty_{2k-1})} \cap
  L^-$, or on both can be assumed to be valued into
  $S^{2k-3}/\mathbb{PC}^{k-2}$.}\\


\noindent{\it Proof of Lemma 3.}
  \medskip
  Assume that $z^\infty_{2k-1}$ has eg at least one large positive $v$-jump. We
  then claim that, for a suitable pseudo-gradient, $\partial
  z^\infty_{2k-1} \cap L^-$ is empty for a large enough index.

  Indeed, let us assume that $z^\infty_{2k-1}$
  dominates $z^\infty_{2k-2}$, of index $(2k-2)$ and that
  $W_u(z^\infty_{2k-2})$ is entirely contained in $L^-$. It follows that $z^\infty_{2k-2}$ has an $H^1_0$-index [3], p7, see also p77, equal
  to zero. For $k$ large, by [3], Lemma 11, p96, $z_{2k-2}^\infty$ must have, after $C^2$-perturbation of the contact form, some
  characteristic (see eg [3], p101)
  $\xi$-pieces. We may
  assume that no decreasing pseudo-gradient may be created at
  $z^\infty_{2k-2}$ with the introduction of a small negative
  $v$-jump anywhere, so that all the characteristic $\xi$-pieces of
  $z^\infty_{2k-2}$ have decreasing normals [4] with the positive
  orientation along $+v$.

  We then introduce a small negative $v$-jump as a companion to
  the now small positive $v$-jump inherited from $z^\infty_{2k-1}$.
  Together, these small negative and positive $v$-jumps can travel
  across the large negative $v$-jumps of $z^\infty_{2k-2}$, until
  the small positive $v$-jump reaches the position of a decreasing
  normal along a characteristic $\xi$-piece of $z^\infty_{2k-2}$ so
  that the flow-line continues past $z^\infty_{2k-2}$, not in $L^-$.
  This characteristic $\xi$-piece must exist for $k$ large enough
  after adjustment of $v$-rotation along $z_{2k-2}^\infty$, see [3],
  Lemma 11, p96. The claim follows and extends with the introduction of
  additional pairs of tiny positive and negative $\pm v$-jumps (this
  does not affect $L^+$ and this does not affect $L^-$)to all
  flow-lines from $z^\infty_{2k-1}$ to $z_{2k-2^\infty}$.  This corresponds to a modification of the pseudo-gradient
flow, from $z^\infty_{2k-1}$, as it reaches $z^\infty_{2k-2}$.

We then claim that $H= \underset {z^\infty_{2k-2} \in
\partial z^\infty_{2k-1}} {\mathrm {\cup}} \overline {W_u(z^\infty_{2k-2}) \cap W_s(L^-\smallsetminus \tilde
J^{-1}(0,\epsilon))}$ can be deformed on a CW-complex of top
dimension $(2k-3)$. This follows from the fact that, above the level
$\epsilon$, the only critical point (at infinity) of $\tilde J$ of
index $1$ is $x^{1, \infty}_-$ and all its other critical points (at
infinity) are of index $2$ or more. Since the $z^\infty_{2k-2}$s are
of index $(2k-2)$, we can use the reverse flow to the decreasing
pseudo-gradient on $H$ and deform it to a CW-complex of dimension
$(2k-3)$.

It follows that we can assume that the classifying map for the
$S^1$-action on $H$ is valued in $\mathbb{PC}^{k-2}$. The claim of
Lemma 3 is established since the additional pieces that we can find
in  $\overline {W_u(z^\infty_{2k-1})} \cap
  L^-$, outside of $H$, are of top dimension $(2k-3)$.

The above proof requires some further work if $z^\infty_{2k-2}$ is
in $\partial^\infty c_{2k-1}$: Indeed, let us consider, for a given
$z^\infty_{2k-2}$, $\overline{W_u(z^\infty_{2k-2}) \cap W_s(x^{1,
\infty}_-)}$. This latter set divides the set $F$ of flow-lines
originating at $z^\infty_{2k-2}$ and abutting to
$J_0^{-1}(\epsilon)$ from the set of flow-lines originating at
$z^\infty_{2k-2}$ and abutting to $B_0 \cap L^-$.

 When $z^\infty_{2k-2}$ is part of $\partial^\infty c_{2k-1}$, the
classifying map is given by the map $"b"$ of [5] on $F\smallsetminus
{z^\infty_{2k-2}}$. The above argument is insensitive to this and we
therefore need in this case a slightly more involved argument,
understanding better the set $H$ introduced above, see below. $\blacksquare$\\

 \noindent {\bf 8. Isotopy of decreasing
  pseudo-gradients}\\

    \noindent  We recall that we have split the critical points at infinity composing $h_{2k-1, \infty}$ into two
   groups, the first group have some large positive $v$-jump,
   whereas the second group has only large negative $-v$-jumps.
  Observe that if $z_{2k-1}^\infty$ has some positive large $v$-jump and
  $k$ is large, then $\overline {W_u(z_{2k-1}^\infty)}\cap L^-$ has,
  according to Lemma 3 above, a classifying map valued into
  $S^{2k-3}/\mathbb{PC}^{k-2}$,
  whereas we can take $L^-$ to be $L^+$ in the above statement if
  $z_{2k-1}^\infty$ has some negative large $-v$-jump.

  Thus, applying Lemma 3 above to our set of specific critical points at infinity, the $y_{2k-1,j}^
  \infty$ of the first group are such that $\overline {W_u(y_{2k-1,j}^\infty)}\cap L^-$
  has a classifying map taking its values
   in  $S^{2k-3}/\mathbb{PC}^{k-2}$; whereas for the second group the second group, it is the classifying
map
  for $\overline {W_u(z_{2k-1,j}^\infty)}\cap L^+$ that is valued into a
  low $S^{2k-3}/\mathbb{PC}^{k-2}$. Completing tangencies as stated above, we may
  assume that the second group is reduced to a single
  $z^{\infty,-}_{2k-1}$.

  We then claim that we can complete, under
  our basic assumption-which we use here in an essential way-an
  isotopy of the decreasing pseudo-gradient which leaves the flow-lines in $L^+$ and
  $L^-$ undisturbed and such that the following claims hold true:

   \medskip
\noindent{\bf Lemma 4} \,  {\it $W_u(z_{2k-1}^{\infty, -}) \cap
D_1^+$ and $W_u(z_{2k-1}^{\infty,
  -}) \cap D_1^-$ are connected.}\\


  \medskip
\noindent {\bf Lemma 5} {\it (i) $W_u(y_{2k-1, j}^{\infty}) \cap
D_1^+$ is
  connected.

  (ii)The classifying map on $\overline {W_u(y_{2k-1,
  j}^{\infty})
  \cap L^-}$ may be assumed to be valued into $S^{2k-3}
  /\mathbb{PC}^{k-2}$.}\\

\noindent {\it Proof of Lemma 4.} \,The arguments for this proof are
strongly inspired from J.Milnor's proof of the h-cobordism theorem,
see Theorem 6.4, p70 of [11].

 We recall that we make the basic
assumption that there are no critical points (at
  infinity) of index $1$, $\tilde x_{\pm}^{1(\infty)}$ connecting
  curves of $J^{-1}(\epsilon) \cap L^+$ and $J^{-1}(\epsilon) \cap
  L^-$. Under our basic assumption, after completing tangencies that
  leaves $L^-$ invariant, we may assume that there is only one
  critical point at infinity of index $1$ $x_-^{1 \infty}$
  connecting $L^-$ and the "small" contractible curves (as above) of $C_\beta$ as
  well as one critical point at infinity of index $1$ $x_+^{1
  \infty}$ connecting $L^+$ and the "small" contractible curves of
  $C_\beta$ (as above), whereas there is no critical point (at infinity) $x^{1
  (\infty)}$ connecting $L^-$ and $L^+$.

  We consider a/the critical point at infinity $z_{2k-1}^{\infty,
  -}$ from $h_{2k-1, \infty}$, as above,
  such that its larger $\pm v$-jumps are all negative and we
  consider a level $c$ just below $J(z_{2k-1}^{\infty, -})$.

 $W_s(L^-) \cap J^{-1}(c)$ is an open
  connected set with a boundary $(\partial W_s(L^-)) \cap J^{-1}(c)$
  that is connected in its top dimension.

  We claim that, for such a critical point at infinity $z_{2k-1}^{\infty, -}$ with only negative large $(-v)$-jumps, we can arrange so that, for
  each $c \lneq J(z_{2k-1}^{\infty, -})$, $c$ close to $J(z_{2k-1}^{\infty,
  -})$, $(W_u(z_{2k-1}^{\infty, -})\smallsetminus W_s(L^-)) \cap
  J^{-1}(c)$ is connected. It suffices for this conclusion that each connected component of
  $W_s(L^-) \cap W_u(z_{2k-1}^{\infty, -}) \cap J^{-1}(c)$ has a
  connected boundary.

  If a connected component, an open set in $W_s(L^-) \cap W_u(z_{2k-1}^{\infty, -}) \cap
  J^{-1}(c)$, has a boundary made of two or more distinct connected
  components $C_1$ and $C_2$, we need to modify the flow, keeping
  the curves of $L^-$ in $L^-$, so that, for this modified flow,
  $C_1$ and $C_2$ are changed and define now the same connected
  component of the boundary of the intersection set.

  The level $c$ is very close to $J(z_{2k-1}^{\infty, -})$ and
  therefore, $C_1$ and $C_2$ may be assumed to be contained in
  $W_s(x_-^{1,
  \infty})$, where $x_-^{1, \infty}$ is the only critical point at
  infinity of index $1$ connecting $L^-$ and the small contractible
  curves of $C_\beta$. We now connect $C_1$ and $C_2$ with two paths $p_1$ and $p_2$, one
  in  $W_s(L^-) \cap W_u(z_{2k-1}^{\infty, -}) \cap
  J^{-1}(c)$, the other one in $W_s(x_-^{1,
  \infty})\cap J^{-1}(c)$. Assuming that $M^3$ is $S^3$, or assuming that $J^{-1}(c)$ is
  connected and simply connected, we may find a surface $\Sigma$ in
  $J^{-1}(c)$ connecting $p_1$ and $p_2$.

  The curves of $W_u(z_{2k-1}^{\infty, -})\cap J^{-1}(c)$ that are
  in $L^-$ define, for $c$ close to $z_{2k-1}^{\infty, -}$, an open
  ball with a connected boundary. We may define our pseudo-gradient
  so that a small open neighborhood of this closed ball flows into
  $L^-$. We may then assume that $p_1$ does not intersect the closure
  of this open ball. In addition, $\Sigma$ may be assumed to be embedded in $J^{-1}(c)$, using
  general position. Again, using general position, $\Sigma$ defines
  the trace of a deformation along which $p_1$ of $W_s(L^-) \cap W_u(z_{2k-1}^{\infty, -}) \cap
  J^{-1}(c)$ is brought on $p_2$. $p_1$ and $p_2$ do not intersect $L^-$. After perturbation, $\Sigma$ also may be assumed
  not to intersect $L^-$: indeed, $\Sigma$ may be assumed to be in some $\Gamma_{2m}$ for $m$ large.
  We may add to the curves of $\Sigma$ $4m$ tiny positive $v$-jumps that are brought to be zero when reaching $p_1$
  and $p_2$. The curves are not in $J^{-1}(c)$ anymore, but they are at a very close level and we can
  flow them back to this level, since none of the curves of $\sigma$
  was critical to begin with. Then, $\Sigma$ does not intersect
  $L^-$.
  This simple deformation can now be "opened
  up" and transformed into an isotopy of decreasing pseudo-gradient.
  At the time $1$ of the deformation, the two modified $C_1$ and
  $C_2$ are now connected, whereas the evolution of the curves of $L^-$ is not
  disturbed.

\medskip


\begin{figure}[!htb]
 \centering
 \includegraphics[scale=.35]{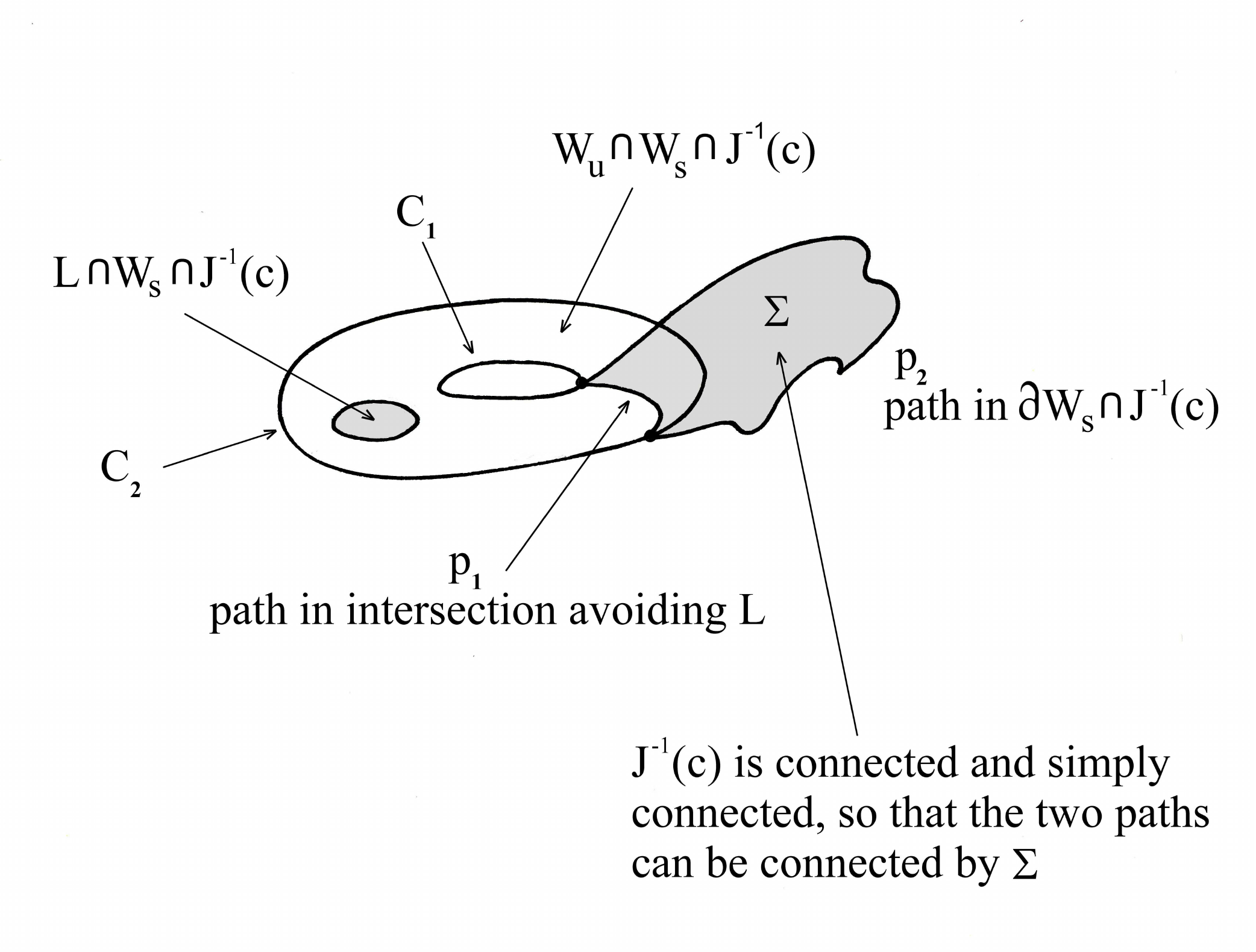}
 \end{figure}

  A similar construction/deformation may be built in order to
  connect  all the various components of $W_u(z_{2k-1}^{\infty, -})$
  going into $L^-$. Once these modifications are performed, we can complete tangencies
  between various $z_{2k-1}^{\infty, -}$s. As long as the tangencies
  occur as described in the figure above, without involving
  flow-lines abutting in $L^-$, the recomposition of the unstable
  manifolds obeys the rule that each connected component of curves
  attracted by $L^-$ has a connected boundary, so that the
  complement of $W_s(L^-)$ in $W_u(z_{2k-1}^{\infty, -})$ (after
  tangencies) is connected (in its top dimension).

The conclusion
  follows for the first claim of Lemma 4. The proof of the second claim follows from the same argument, slightly modified. $\blacksquare$

\noindent {\it Proof of Lemma 5.}\, The only statement that requires
additional
  proof is the claim about the classifying map. The addition of the
  various $\Sigma$s built as above does not change the
  Fadell-Rabinowitz index since these are contractible pieces and
  they may be assumed not to dominate any critical point above
  $D_1^-$ (after re-parametrization, see above and J. Milnor [11]).
  Then, after "opening up $\Sigma$" as above, we find that $\overline {W_u(y_{2k-1,
  j}^{\infty})
  \cap L^-}$ is contained in a set having a classifying map valued
  into $S^{2k-3}
  /\mathbb{PC}^{k-2}$ as claimed. $\blacksquare$\\


  The arguments collapse if $\partial W_s(L^-)$ is not connected.\\

 \noindent { \bf  9. The extension of Lemma 1 to $\overline
  {W_u(h_{2k-1, \infty})}\cap \tilde
  J^{-1}[\epsilon, \infty)$}\\


\noindent {\bf Proposition 1} {\it  (i)  Lemma $1$ extends to
$\overline
  {W_u(h_{2k-1, \infty})}\cap \tilde J^{-1}[\epsilon, \infty)$. The
  classifying map after deformation is valued into $S^{2k-3},
  \mathbb{PC}^{k-2}$.

  (ii)  Along this deformation, the classifying map restricted to $\overline
  {W_u(h_{2k-1, \infty})}\cap (\tilde J^{-1}(\epsilon)\cup W_u(x^{1, \infty}_-))=\overline
  {W_u(h_{2k-1, \infty})}\cap B_0$ is valued into $(\mathbb{PC}^{k-1}
  \times [0,1] \cup \mathbb{PC}^{k-2} \times [-1, 1] \cup
  \mathbb{PC}^{k-1} \times \{-1\})$.}\\


\noindent {\it Proof of Proposition 1.}\\

  \noindent  {\it Extending Lemma 3 to $\partial^\infty c_{2k-1}$, with the
  $"b"$ pre-assigned value [5] of the classifying map when the
  $v$-component of the tangent vector to the curves has at least one
  sign-change.}\\

  In a first step, we extend Lemma 3 and we prove that, if $y_{2k-2}^\infty$ is in $\partial^\infty c_{2k-1} \cap \partial y^\infty_{2k-1, j}$,
  then the classifying
  map $"b"$ of [5] can be {\bf extended} to $\overline{W_u(y_{2k-1,j}^\infty) \cap L^-}$
  with values in $S^{2k-3}/\mathbb{PC}^{k-2}$ on this latter set.

  We need for this a few preliminary definitions, Lemmas etc. We
  start with:\\


 \noindent {\bf Definition 1}  \,{Let $(\partial^\infty c_{2k-1})_-$ be the
critical points at infinity in $\partial^{\infty} c_{2k-1}$ having
all their large $v$-jumps oriented along $-v$ and having a non-zero
$H^1_0$-index.}

\medskip

\noindent {\bf Requirements on decreasing flow-lines.}
   {\it We are requiring that our decreasing pseudo-gradient leaves the sets
$L^+$ and $L^-$ invariant (respectively) and that it never increases
the number of zeros of the $v$-component $b$ of the curves under
decreasing deformation, this {\bf solely} for closure of the set of
flow-lines originating at any periodic orbit of index $(2k-1)$.}
\medskip

Therefore, starting from $y^\infty_{2k-1, j}$ as above, which has at
least one large positive $v$-jump, and reaching to a critical point
at infinity of $(\partial^{\infty} c_{2k-1})_-$, we find curves that
have a mixture of positive and of negative steady $\pm v$-jumps. On
such curves, we can add additional negative or positive $\pm
v$-jumps as we please, we are not bound by any requirement since the
flow-line is not originating at a periodic orbit of index $(2k-1)$.
\medskip
We then claim: \

\noindent{\bf Lemma 6}  \,\,{\it Any critical point at infinity in
$(\partial^{\infty} c_{2k-1})_-\cap \partial z^{\infty}_{2k-1,j}$
has no characteristic $\xi$-piece. After a $C^2$-bounded,
$C^1$-small perturbation of the contact form $\alpha$ in the
vicinity of this critical point at infinity, we may assume that the
maximal number of sign-changes for $b$ on its unstable manifold is
$(2k-4)$.}\\


\noindent {\bf Remark 1}\,\, Lemma 6 is not absolutely required in
our proof
of Theorem 1, but it is a convenient result.\\


\noindent {\it Proof of Lemma 6.} Following our requirements and
observation above, this critical point at infinity cannot have any
characteristic $\xi$-piece, since we are then free, on flow-lines
out of $z_{2k-1,j}^\infty$ and reaching this critical point at
infinity, to introduce a decreasing normal [4] along this
characteristic $\xi$-piece and bypass this critical point at
infinity. We may assume that it has some non-zero $H^1_0$-index for
$k$ large enough. Indeed, otherwise, we can use Lemma 3 and
Proposition 15 of [3]. There is enough $v$-rotation on the various
$\xi$-pieces and we can transport it in a given $\xi$-piece, thereby
creating a non-zero $H^1_0$-index on this $\xi$-piece.

Since $c_{2k-1}$ dominates this critical point at infinity, the
maximal number of zeros on its unstable manifold is $(2k-2)$ at
most. Since it has a non-zero $H^1_0$-index, we can use Lemma 3 of
[3] and modify at least by $2$ this maximal number of zeros. The
claim follows.$\blacksquare$\\

 \noindent {\bf Lemma 7} {\it  $x^{1, \infty}_-$ may be assumed to have at least
one large positive and one large negative $v$-jump.}\\


\noindent {\it Proof of Lemma 7. } $x^{1, \infty}_-$ introduces a
genuine difference of topology in the level sets of the functional
$J$. It cannot therefore have a characteristic $\xi$-piece. Would it
have eg only negative large $v$-jumps, then its $H^1_0$-index cannot
be zero: $x^{1, \infty}_-$ connects $J^{-1}_0(\epsilon)$ and
$J^{-1}(\epsilon)\cap L^-$ and this cannot be achieved with a Morse
index totally at infinity.

Since the $H^1_0$-index of $x^{1, \infty}_-$ is non-zero, we can
modify it using again Lemma 3 of [3]. It cannot become $2$, this
would be too high. Thus, it has to become zero; the index of
$x^{1,\infty}_-$ is totally at infinity and this is a contradiction
as
pointed out above. $\blacksquare$\\

It follows that there exists a neighborhood of $W_u(x^{1, \infty}_-)
\cap J^{-1}([\epsilon, \infty))$ where the classifying map for the
$S^1$-action may be assumed to be given by the map $"b"$ of [5],
since the $v$-component of $\dot x$ has at least two zeros.


\medskip
 {\it The classifying map on $\overline{\underset {z^{2, \infty}_l}
{\mathrm {\cup}} W_u(\partial^{\infty} c_{2k-1})_-\cap
\partial y^{\infty}_{2k-1,j}\cap W_s(z^{2, \infty}_l)}$  and nearby}
\medskip


 Thus,
 the classifying map is given on part of
$\overline{W_u(\partial^{\infty} c_{2k-1})_-\cap \partial
y^{\infty}_{2k-1,j})\cap W_s(x^{1, \infty}_-)}$ and there is now the
need to extend this map to a set that retracts by deformation on

 $\overline{\underset {z^{2, \infty}_l}
{\mathrm {\cup}} W_u(\partial^{\infty} c_{2k-1})_-\cap
\partial y^{\infty}_{2k-1,j}\cap W_s(z^{2, \infty}_l)}$ where the
$z^{2, \infty}_l $ are all the critical points at infinity of index
$2$ dominating $x^{1, \infty}_-$. This is a stratified set of top
dimension $(2k-4)$. Its classifying map may be assumed, by general
position, to be valued into $\mathbb{PC}^{k-2}$. A homotopy of this
classifying map may also be assumed, using the same general position
argument, to be valued into $\mathbb{PC}^{k-2}$.


The critical points at infinity of this stratified set are of two
types: there are those which contain a sign-change in their large
$\pm v$-jumps. The map "b" of [5] is well-defined on a full
neighborhood of these critical points at infinity.

Then, there are those having all negative large $(-v)$-jumps. Their
$H^1_0$-index cannot be zero since they dominate $x^{1, \infty}_-$
which has a sign-change in its large $\pm v$-jumps. We define in a
neighborhood of these critical points at infinity a $"b"$-map which
is slightly different from the map $"b"$ of [5]: there is a
connected region, diffeomorphic to a cone, in the unstable manifold
of such a critical point at infinity made of curves such all
possible $\pm v$-jumps are non-zero and negative. On the boundary of
this region, some of these negative $v$-jumps are zero. All of these
correspond to $H^1_0$-directions near the dominating critical point
at infinity.

Along this boundary, turning one of the zero $v$-jumps corresponding
to $H^1_0$-index directions into positive tiny $v$-jumps defines a
convex entering set of normal directions into the curves of the
unstable manifold where $b$ changes sign. Furthermore, if this
critical point at infinity dominates another critical point at
infinity of the same family with a non-zero $H^1_0$-index, then,
since all $\pm v$-jumps that are non-zero on this boundary are
negatively oriented, we derive that this $H^1_0$-position must have
existed above, in the dominating critical point at infinity and must
have survived all along the flow-lines connecting these two critical
points at infinity of the same family. It follows that the set of
entering normals is well-defined. Since the regions where all
possible $\pm v$-jumps are negative cannot dominate $x^{1,\infty}$,
we find that we can use this set of entering normals and define the
map $"b"$ all over our stratified set, except for the periodic
orbits. Observe that, if on some flow-lines originating at one of
the critical points at infinity as above, with all large negatively
oriented $\pm v$-jumps, there is a positive $v$-jump due to the use
of an $H^1_0$-direction and that this positive $v$-jump cancels with
a negative $-v$-jump as we approach a lower critical points of the
same family, then the map $"b"$ of [5] is defined on the flow-lines,
originating and ending critical points at infinity excluded. Using
Lemma 6 above, it can be glued with the map $"b"$ as defined above,
with values into $S^{2k-3}/\mathbb{PC}^{k-2}$. Observe in addition
that if, starting from $y^\infty_{2k-1, j}$, we end up at a critical
point at infinity of $\partial^{\infty} c_{2k-1}$ with all its $\pm
v$-jumps oriented along $+v$, then the map $"b"$ of [5] will be
defined in the vicinity of the flow-lines starting at this critical
point at infinity and ending into $L^-$, with at least two zeros and
at most $(2k-4)$ zeros and we can again glue this map with the other
map $"b"$ as defined above, with a resulting map valued in
$S^{2k-3}/\mathbb{PC}^{k-2}$. We could use a weaker statement than
the statement of Lemma 6, with $(2k-2)$ zeros in lieu of $(2k-4)$.

The periodic orbits are of top index $(2k-3)$, with a maximal number
of zeros of $b$ on their unstable manifold equal to $(2k-4)$. In
order to define the map $b$, we need $b$ to have at least two zeros.
$b$ is identically zero at the periodic orbit, but we can perturb
the unstable manifold so that $b$ is non-zero at the top perturbed
critical point and has $(2k-2)$ zeros, with a maximal number of
zeros for $b$ on this perturbed unstable manifold equal to $(2k-2)$
near the top, $(2k-4)$ below; this, if the periodic orbit is of
index $(2k-3)$; $(2k-4)$ otherwise, in lieu of $(2k-2)$.
 The
flow-lines that dominate $x^{1, \infty}$ in this unstable manifold
must be such that $b$ has at least two zeros on their curves. There
could be other periodic orbits/critical points at infinity in their
closure, for which we proceed as above.

The resulting map $"b"$ extends to this stratified set, valued into
$\mathbb{PC}^{k-2}$.
\medskip

\noindent {\it Resolving the multiplicity of $\overline{\underset
{z^{2, \infty}_l} {\mathrm{\cup}} W_u((\partial^\infty
c_{2k-1})_-\cap
\partial y_{2k-1,j}^\infty)\cap W_s(z_l^{2, \infty})}$ at the
critical points (at infinity) that it contains.}
\medskip


We now resolve the "multiplicity" of this stratified set of
decreasing flow-lines at each critical point (at infinity), thereby
creating a stratified set $T_{2k-4}$, which is a section to the
decreasing flow abutting into $L^-$.

Indeed, the original set is a closed invariant set of decreasing
flow-lines. Far away from the critical points (at infinity), it can
be perturbed into a section to a decreasing flow abutting into
$L^-$. Close to the critical points (at infinity), we find possibly
several "leaves" for this stratified set, intersecting at the
critical point (at infinity). The "leaves" define components, some
of them abutting to $L^-$, the other ones to eg
$J_0^{-1}(\epsilon)$. We can resolve them also into sections to a
decreasing flow.
\medskip

\begin{figure}[!htb]
 \centering
 \includegraphics[scale=.3]{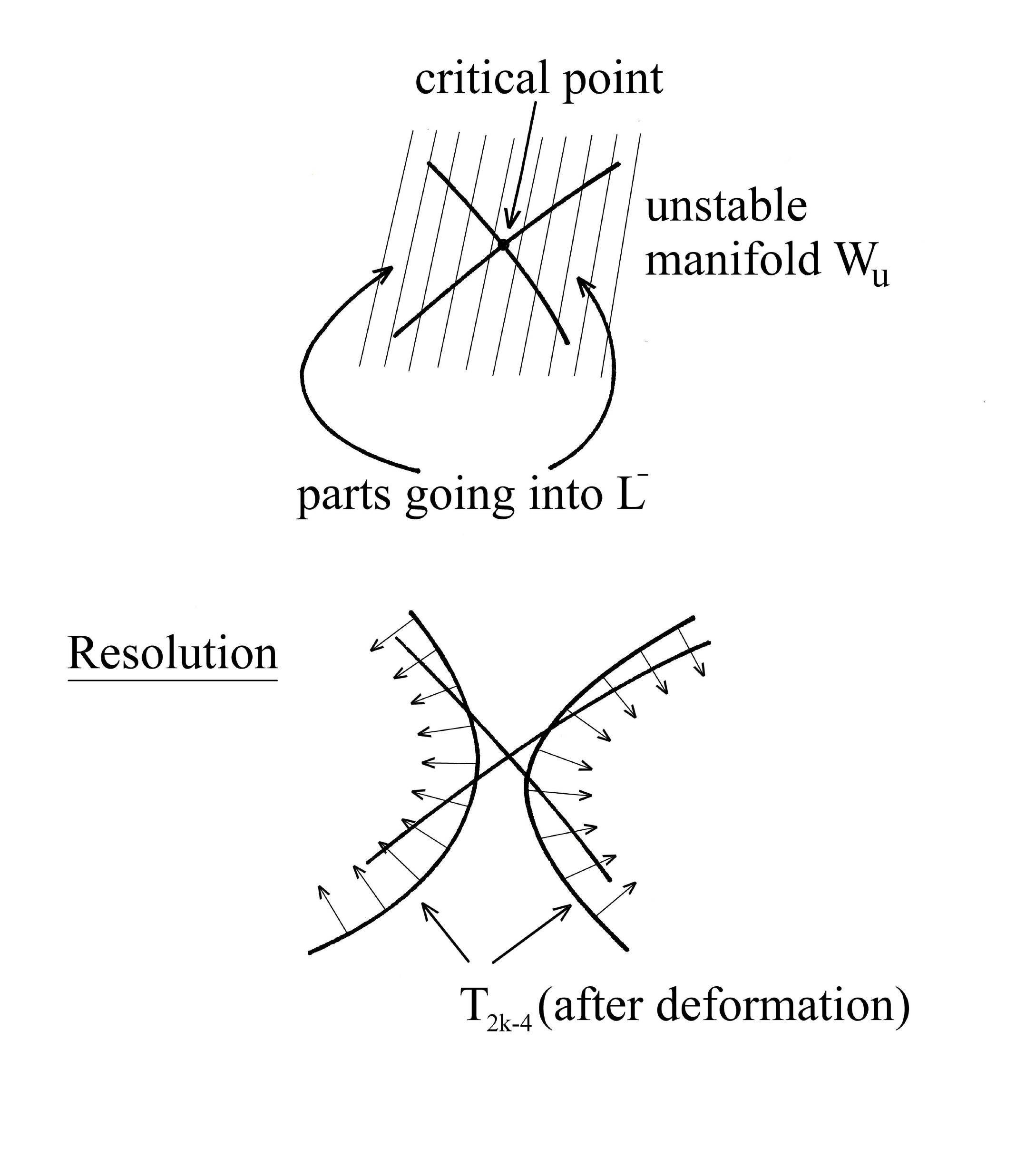}
 \end{figure}


On $T_{2k-4}$, two classifying maps are now defined: the map $"b"$
as above and the map $\Psi$, mapping $T_{2k-4}$ into its limit set
at infinity $L^-_\infty$ and from there, to $\mathbb{PC}^\infty$.
$L^-_\infty$ is, after deformation, of top dimension $(2k-4)$, so
that $\Psi$ may be assumed to be valued into $\mathbb{PC}^{k-2}$
(top dimension $(2k-3)$ would lead to the same conclusion).

The homotopy between these two maps $"b"$ and $\Psi$, restricted to
$T_{2k-4}$ may be assumed to be valued into $\mathbb{PC}^{k-2}$ as
well.

We now conclude the argument. The figures of reference are as
follows:
\medskip

\begin{figure}[!htb]
 \centering
 \includegraphics[scale=.4]{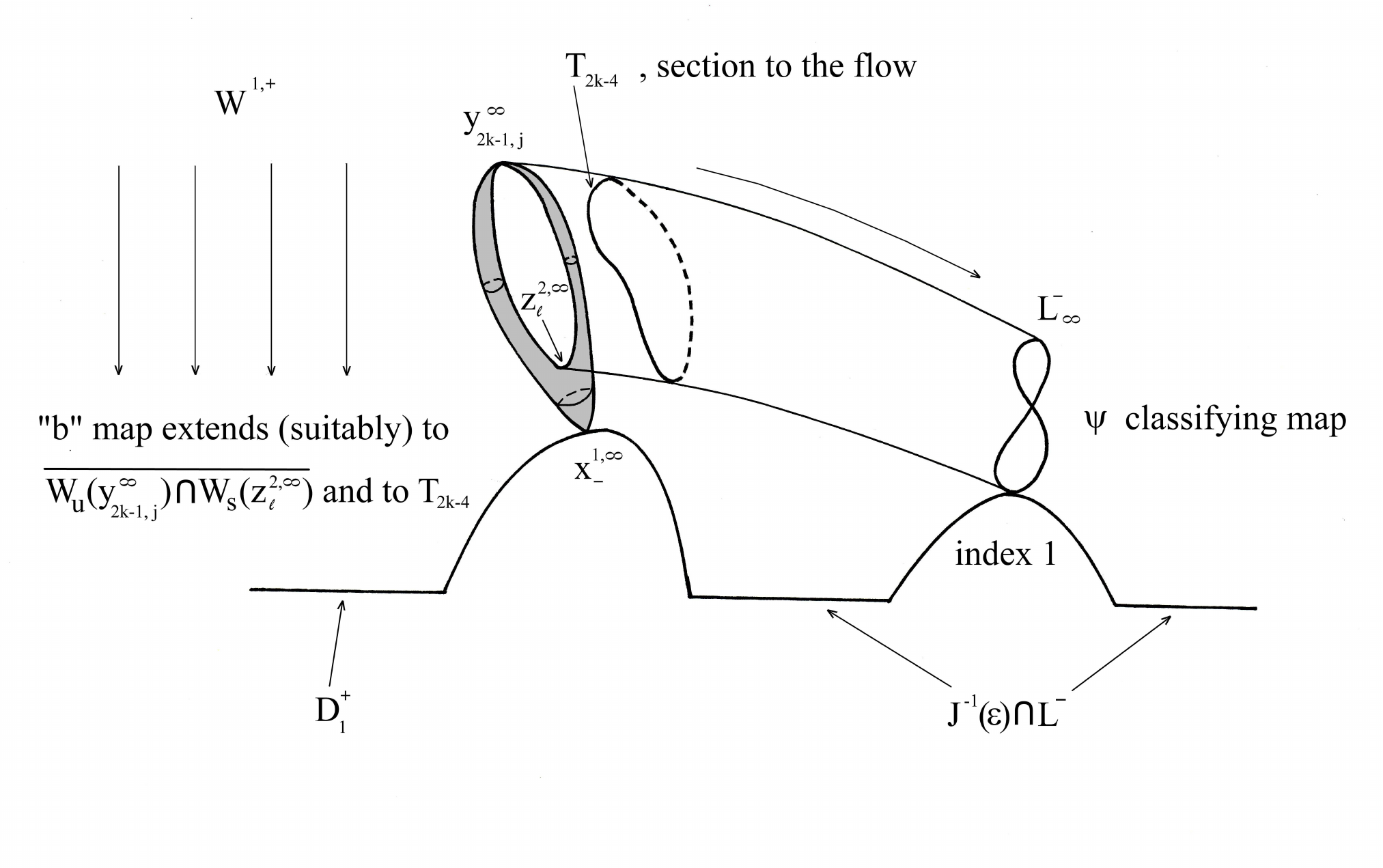}
 \end{figure}
\medskip

\begin{figure}[!htb]
 \centering
 \includegraphics[scale=.4]{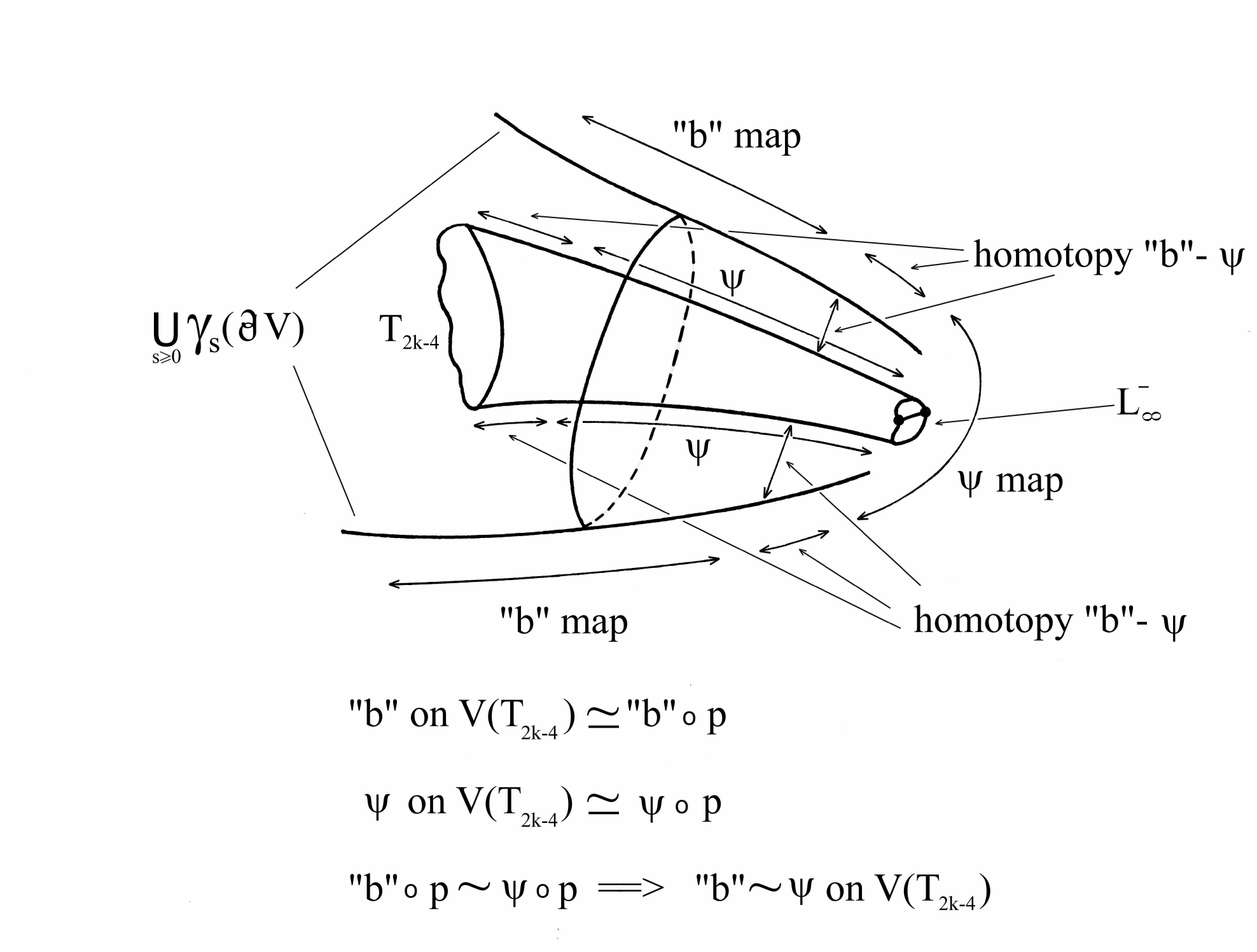}
 \end{figure}


Taking a small neighborhood $V$ of $T_{2k-4}$ in section to the
flow, we may flow it using the decreasing flow $\gamma_s$.
$\underset { s \geq 0} {\mathrm {\cup}} \gamma_s(V)$ has $\underset
{ s \geq 0} {\mathrm {\cup}} \gamma_s(\partial V)$ as a boundary.
There is a projection $p:\underset { s \geq 0} {\mathrm {\cup}}
\gamma_s(V) \longrightarrow \underset { s \geq 0} {\mathrm {\cup}}
\gamma_s(T_{2k-4})$ and the map $"b"$ and the map $\Psi$ which were
defined above on $T_{2k-4}$ thereby extend, using $p$, to $\underset
{ s \geq 0} {\mathrm {\cup}} \gamma_s(V)$. Observe that the standard
map $"b"$ of [5] is homotopic to the map $"b"$ defined above, the
homotopy being valued into $\mathbb{PC}^{k-2}$ and observe that, on
$\underset { s \geq 0} {\mathrm {\cup}} \gamma_s(\partial V)$, this
map maybe viewed after deformation as constant on each $\underset {s
\geq 0}{\mathrm {\cup}} \gamma_s(z)$, equal to its value at $p(z)$,
since on each flow-line the changes of sign of the function $b$ can
be recorded unchanged as the time-parameter $s$ increases; it is
only the sizes of the function $b$ and its shapes, on the various
intervals between zeros that we can track, which change.

The $\omega$-limit set of each $\underset { s \geq 0} {\mathrm
{\cup}} \gamma_s(V)$ and $\underset { s \geq 0} {\mathrm {\cup}}
\gamma_s(\partial V)$ is the same: it is $L^-_\infty$, with its map
$\Psi$. Since the map $"b"$ and the map $\Psi$ are homotopic when
restricted to $T_{2k-4}$, with a homotopy valued into
$\mathbb{PC}^{k-2}$ and since their values on $\underset { s \geq 0}
{\mathrm {\cup}} \gamma_s(V)$ are derived with the use of $p$, they
are homotopic as maps defined on this larger set, with the same
target value set $\mathbb{PC}^{k-2}$.

As we reach to $L^-_\infty$, starting with $\partial V$ and flowing
down, we may gradually use this homotopy and insert the map $\Psi$,
so that the classifying map takes the well-defined value $\Psi$ on
$L^-_\infty$. Going deeper into $\underset { s \geq 0} {\mathrm
{\cup}} \gamma_s(V)$, we use more and more the map $\Psi$ on the
flow-lines. When we reach $T_{2k-4}$, the map is $\Psi$ all along
the decreasing flow-lines. Of course, we have used an interval
$[-\epsilon,0]$ of times $s$ to replace $"b"$ by $\Psi$ as we start
in $T_{2k-4}$.

We have therefore extended the map $"b"$ on $\partial^\infty
c_{2k-1} \cap \partial z^\infty_{2k-1, j}$ to the flow-lines
abutting in $L^-$ and the extension is valued into
$\mathbb{PC}^{k-2}$. Using the fact that the "bottom set" $D_1^+$ is
connected, we may now apply, without perturbing the topological
arguments of section 11, below, the procedure of Lemma 1 above to
the topological boundary of $W_u(y^\infty_{2k-1, j}) \cap
J^{-1}([\epsilon, \infty))$. We find a classifying map valued into
$\mathbb{PC}^{k-2}$ on $\overline {W_u(y^\infty_{2k-1, j}) \cap
J^{-1}([\epsilon, \infty))}$. We will use this later.\\

 {\it Conclusion for the extension of Lemma 3.}\\

 We complete the modifications described in the first part of
  this paper, for all $W_u(y_{2k-1, j}^\infty)$s such that
  $\partial y_{2k-1, j}^\infty \cap L^-$ has a classifying map valued into $\mathbb{PC}^{k-2}$. The modifications do not occur
  on
  flow-lines abutting in $L^-$ then since, by Lemma 3, the classifying map on $\overline{W_u(y_{2k-1, j}^\infty)} \cap
  L^-$, and even on $\overline{\underset {m} {\mathrm {\cup}} W_u( y_{2k-1, m}^\infty )}\cap
  L^-$, may be assumed to be given, valued in $S^{2k-3}, \mathbb{PC}^{k-2}$. These modifications occur on flow-lines abutting in $D_1^+$. We
  know that each $\partial W_u(y^\infty_{2k-1, j})$ is connected. By
  Lemma 5, we know that $W_u(y^\infty_{2k-1, j})\cap D_1^+$ is
  connected and, according to the construction of $D_1^+$, see
  section 4, no critical point (at infinity) of index $1$ dominates
  $D_1^+$, aside from $x^{1, \infty}_-$.

  The arguments for Lemma 1 can then be applied to each of these $W_u(y_{2k-1,
  j}^\infty)$s.

  Once the classifying map is defined on these unstable manifolds in $h_{2k-1, \infty}$, we are left with the $z_{2k-1, j}^\infty$
   of $h_{2k-1, \infty}$ such that their
  large $\pm v$-jumps are along $-v$. We have reduced them to a
  single $z_{2k-1,-}^\infty$, which we denote $z^\infty$ in the
  sequel.\\

{ \it The conclusion for the proof of Proposition 1.}\\

Let now $\overline {W^{1,+}}$ be the closure of the set of
decreasing flow-lines abutting to the "bottom set" $D_1^+$.

Arguing as above, but using $z^{\infty,-}_{2k-1}$ in lieu of
$y^\infty_{2k-1, j}$, we may assume that the classifying map on
$\overline{\partial ^\infty W_u(c_{2k-1}) \cap \partial W_u(
z^\infty_{2k-1,-})\cap W^{1,+}}$ is also valued into
$\mathbb{PC}^{k-2}$: this involves extending as above a variant of
the map $"b"$ of [5] into $L^+$. The reasoning is identical to the
case  for $y^\infty_{2k-1,j}$, only that $L^-$ is now replaced with
$L^+$.

 There is however no global reduction of the classifying map on all of
$\overline {W_u(z^{\infty,-}_{2k-1})}$ as above for $\overline
{W_u(y^\infty_{2k-1, j})}$ since the "bottom set" is not connected
now. The argument is different. It goes as follows:

After our reasoning above, also Lemma 1 and Proposition 1, we know
that the classifying map is valued into $\mathbb{PC}^{k-2}$ on
$\overline {W_u(y^\infty_{2k-1, j})}$, on the trace of
$h_{2k-1}^\infty$ and $c_{2k-1}$ on the bottom set $D^+_1$ and also
on $\overline{\partial ^\infty W_u(c_{2k-1}) \cap \partial W_u(
z^\infty_{2k-1,-})}\cap \overline{ W^{1,+}}$. Since $\partial
y_{2k-1,j}^\infty +\partial z_{2k-1,-}^\infty +\partial ^\infty
c_{2k-1}=0$, we derive from the claims above that the classifying
map on $\overline{\partial W_u(z^\infty_{2k-1,-})\cap W^{1,+}}$ is
valued into $\mathbb{PC}^{k-2}$. Using the proof of Lemma 4 and the
proof of Lemma 5 and the connectedness of $\overline{
W_u(z^\infty_{2k-1,-})\cap
\partial \overline {W^{1,+}}}$, we derive, since this set and
$\overline{\partial W_u(z^\infty_{2k-1,-})\cap \overline{ W^{1,+}}}$
add up to a boundary of top dimension $(2k-2)$, that $\overline{
W_u(z^\infty_{2k-1,-})\cap
\partial \overline {W^{1,+}}}$ has also a classifying map valued into
$\mathbb{PC}^{k-2}$.

  Through our previous modifications, the classifying map
  is given on $(W_u(c_{2k-1}) \cup W_u(h_{2k-1, \infty}\smallsetminus
  z^\infty))\cap D_1^+$, valued into $\mathbb{PC}^{k-2}$.

  This classifying map can be extended to $(W_u(c_{2k-1}) \cup W_u(h_{2k-1, \infty}))\cap
  D_1^+$, valued into $\mathbb{PC}^{k-1}$. By Lemma 2, it is of degree zero. Since
  this map restricted to $(W_u(c_{2k-1}) \cup W_u(h_{2k-1, \infty}\smallsetminus
  z^\infty))\cap D_1^+$ is valued into $\mathbb{PC}^{k-2}$ and since
  $W_u(z^\infty) \cap D_1^+$ is connected, we can modify the
  classifying map relative to this preassigned value on  $(W_u(c_{2k-1}) \cup W_u(h_{2k-1, \infty}\smallsetminus
  z^\infty))\cap D_1^+$ so that it is now valued into
  $\mathbb{PC}^{k-2}$.

  It follows that the {\bf topological} boundary $\partial W_u(z^\infty) \smallsetminus (\partial
  W_u(z^\infty) \cap L^-)$ is of Fadell-Rabinowitz index $(k-2)$ and
  therefore, the {\bf topological} boundary $(\partial
  W_u(z^\infty) \cap L^-)$ is also of Fadell-Rabinowitz index also
  $(k-2)$. By Lemma 4, it is a connected set if we attach to it, without
  increasing its index, boundaries of appropriate neighborhoods (see section 4, above) of unstable manifolds of
  critical points at infinity of index $1$ connecting the various
  components of $J^{-1}(\epsilon) \cap L^-$. These neighborhoods were used in section 4 in order to define the
  appropriate "bottom set" $D_1^-$ in $L^-$, formed essentially of
  $J^{-1}(\epsilon) \cap L^-$ and of these unstable manifolds, glued
  together so that this defines a "level surface" (ie a "bottom set"
  transverse to the flow), see section 4.

  We may therefore assume that, on all of $\partial W_u(z^\infty)$ as well as on the trace of $W_u(z^\infty)$ on $B_0=D_1^+ \cup D_1^-\cup W_u(x^{1, \infty}_-)$,
  the classifying map is given, extending the one previously defined on $\overline{W_u(h_{2k-1, \infty}\smallsetminus
  z^\infty)}$  valued into $S^{2k-3}, \mathbb{PC}^{k-2}$.

  Using the arguments of Lemma 1, this map can now be extended to
  $W_u(z^\infty)$, so that the modifications of Lemma 1 have
  now been completed on all of $\overline{W_u(h_{2k-1, \infty})}$, with a trace
  on the bottom set $B_0$ valued into $(\mathbb{PC}^{k-1} \times \{-1\} \cup \mathbb{PC}^{k-1} \times
  [0,1] \cup \mathbb{PC}^{k-2} \times [-1, 1])$.

Summarizing, the scheme of proof of Theorem 1 is as follows,
supported by the following figure:

\medskip

\begin{figure}[!htb]
 \centering
 \includegraphics[scale=.3]{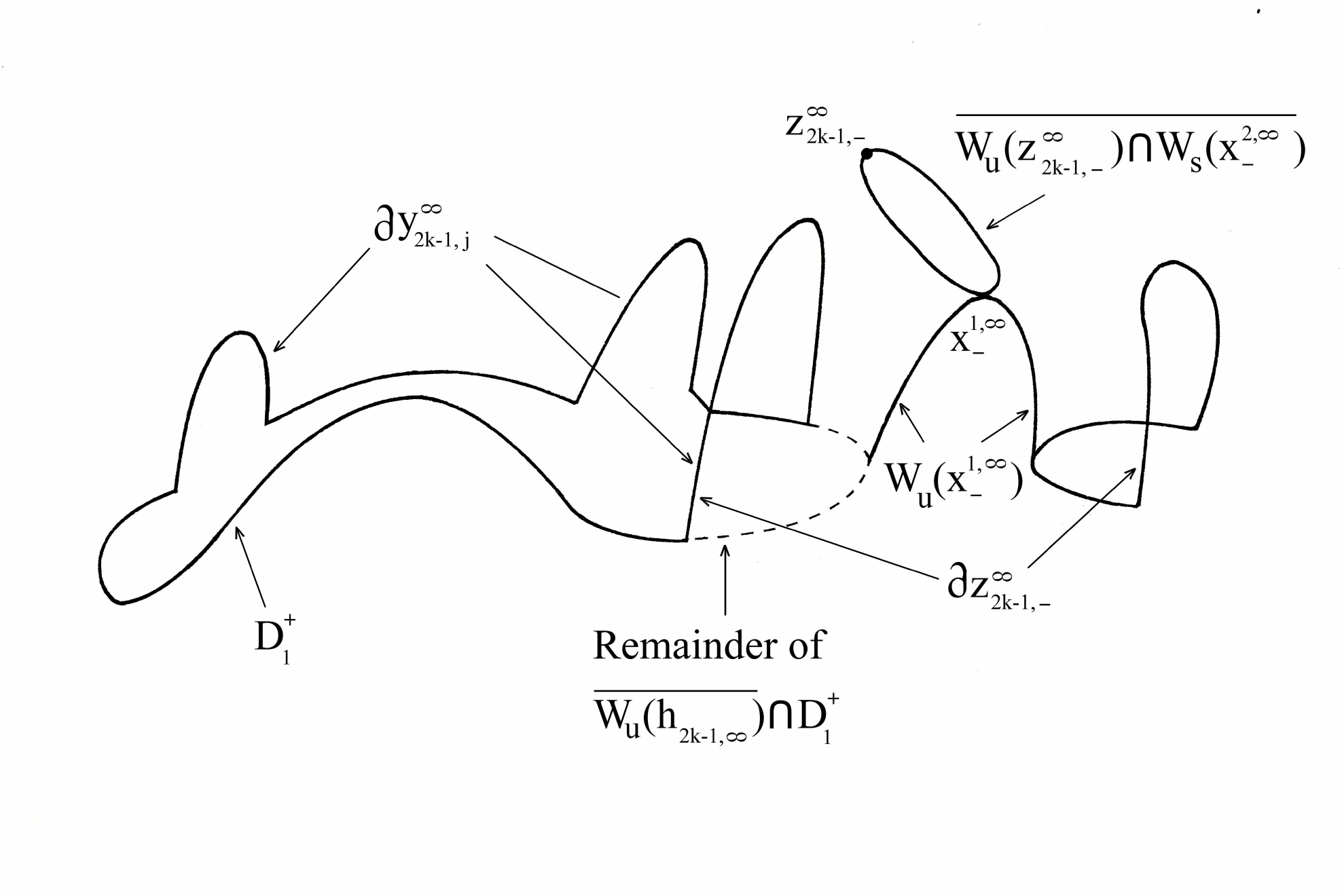}
 \end{figure}


\medskip

.Step1: The classifying map is valued into $\mathbb{PC}^{k-2}$ on
$\partial y_{2k-1,j}^\infty \cup (\overline{W_u(y^\infty_{2k-1,j})
\cap D^+_1)}$ (Lemma 1).

.Step2: The classifying map can be extended to the trace of
$\overline {W_u(h_{2k-1, \infty})}$ on $D_1^+$, valued in
$\mathbb{PC}^{k-2}$. Therefore, the classifying map on
$\overline{W_u(z_{2k-1, -}^\infty)\cap W_s(x^{1, \infty}_-)}$ is
valued in $\mathbb{PC}^{k-1}$, with degree zero.

.Step3: We know that $\overline{\partial z_{2k-1, -}^\infty \cap
W_s(L^-)\cup W_u(z_{2k-1,-}^\infty) \cap D_1^-}$ is of dimension
$(2k-2)$ and connected. From Step 2, we derive that the classifying
map on this set is of degree zero and the conclusion follows.
$\blacksquare$


  \bigskip
\noindent{\bf 10.Multiplicity of domination in dimension $(2k-1)$
and $(2k-2)$, Algebraic Intersection Numbers and
  Flow-lines}

  \medskip

  If a $y_{2k-1, j}^\infty$ appears multiple times in the definition
  of $h_{2k-1, \infty}$, or if $z_{2k-1}^{\infty, -}$ appears a
  number of times, we may resolve this multiplicity and introduce
  several distinct critical points, as many as needed, with very
  close unstable manifolds. The functional is slightly changed and
  its critical points as well, but the arguments are essentially the
  same.

  We need now to resolve the multiplicities of $\overline
  {W_u(h_{2k-1, \infty})}$ at the order $(2k-2)$.
  \medskip

\noindent{\it The case for the $y_{2k-1, j}^\infty$s}.
  \medskip

 Following the technique introduced above, we claim that:\\

\noindent{\bf Lemma 8}\, The decreasing flow can be modified  so
that the algebraic intersection
  numbers $y_{2k-1,j}^\infty-z_{2k-2}^{(\infty)}$ are equal in absolute value to the
  number of actual flow-lines from $y_{2k-1,j}^\infty$ to
  $z_{2k-2}^{(\infty)}$. $L^+$ and $L^-$ remain invariant under this flow.\\


  \noindent{\it Proof of Lemma 8}.
   We need to complete cancellations of flow-lines from a $y_{2k-1,j}^\infty$ to a $z_{2k-2}^{(\infty)}$ with
  opposite intersection numbers $+1$ and $-1$. Between  $y_{2k-1, j}^\infty$
  and
  $z_{2k-2}^{(\infty)}$, for $(2k-2) \geq 2$, we may assume that we
  do not find any critical point (at infinity) of index $1$. After
  re-parametrization of the flow-lines as in [11], Theorem 4.1, pp37-38, there is no
  loss of generality in this assumption. Then, the traces of the unstable manifold
  of $y_{2k-1,j}^{\infty}$ and of the stable manifold of
  $z_{2k-2}^{(\infty)}$ on an intermediate level surface $J^{-1}(c)$ may be assumed to be
  connected. if $M^3=S^3$, we may also assume, without loss of
  generality, that this level surface is simply connected. If $M$ is
  not $S^3$, some more work is required.

  We then join two intersection points with opposite intersection
  numbers in $W_u(y_{2k-1,j}^\infty)\cap J^{-1}(c)$ and in
  $W_s(z_{2k-2}^{(\infty)})\cap J^{-1}(c)$ with two paths $p_1$ and
  $p_2$. We connect $p_1$ and $p_2$ along a surface $\Sigma$, as
  above, in $J^{-1}(c)$. We "slide" as above
  $W_u(z_{2k-2}^{(\infty) })$ along $\Sigma$, modifying it in this
  way. At the end of the process, the cancellation of the two
  intersection points is performed. The argument follows the work of
  J.Milnor (Proof of the h-cobordism theorem) [11], Theorem 6.1, p70. The remaining various boundaries between the
  various critical points at infinity of index $(2k-1)$ can be
  pieced together so that there is no singularity in dimension
  $(2k-2)$ and the argument can proceed.
  \medskip

\begin{figure}[!htb]
 \centering
 \includegraphics[scale=.3]{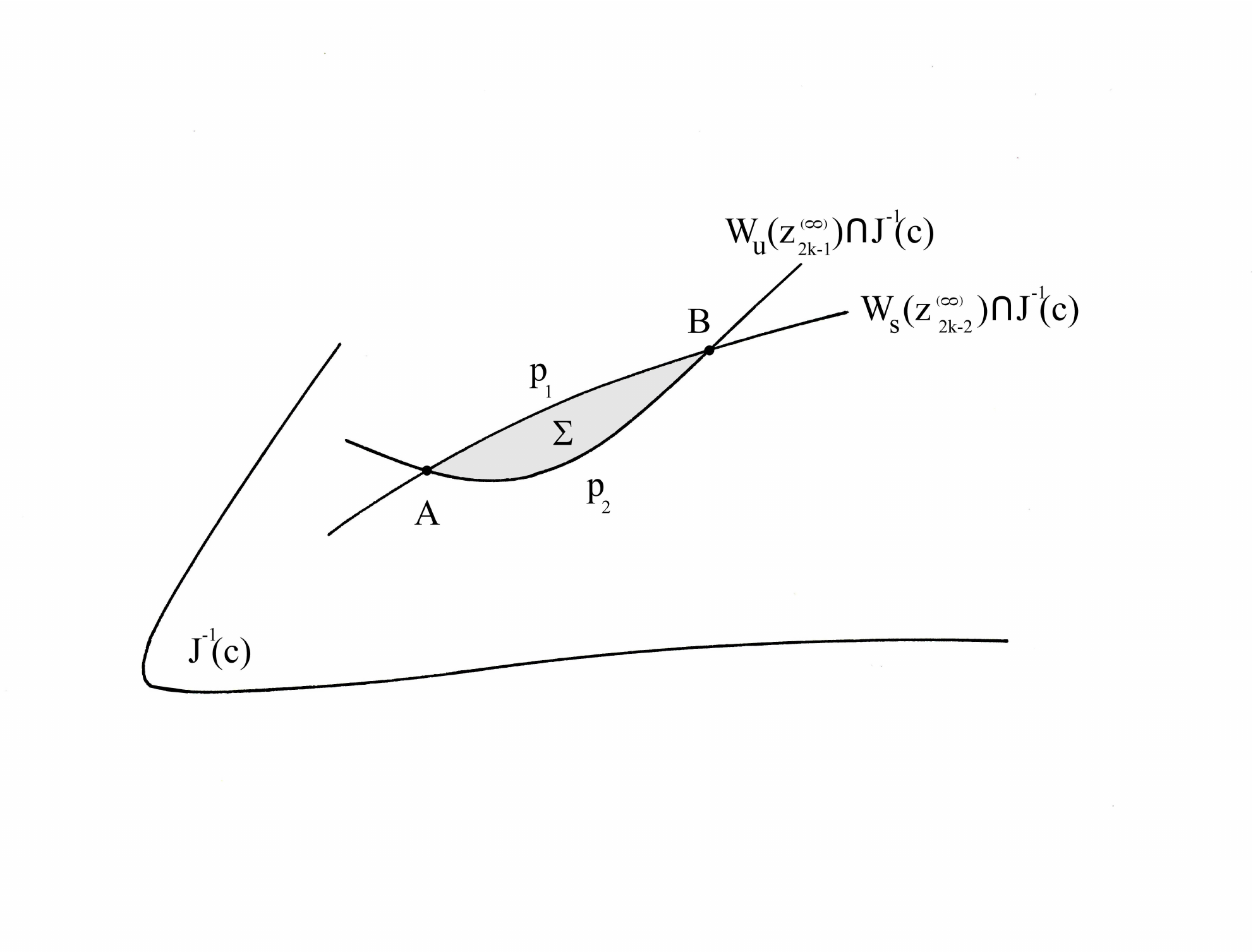}
 \end{figure}


\medskip

  Of course, we need to check that this does not perturb the
  flow-lines in $L^-$. This is quite clear for the $y_{2k-1,
  j}^\infty$s as above. $\blacksquare$

  \medskip

  \noindent{\it The case for $z_{2k-1, -}^\infty$}
  \medskip

  For $z_{2k-1, -}^\infty$, some
  additional care is required. However, we can then modify the
  argument here: if $z_{2k-1, -}^\infty$ dominates a critical point
  at infinity of $L^-$ of index $(2k-2)$ with an algebraic number of
  intersection equal to $0$ with two flow-lines of opposite
  intersection numbers $+1$ and $-1$, we can introduce an additional
  critical point of index $(2k-1)$ and resolve with the help of this
  additional critical point this multiple domination into simple
  dominations of distinct critical points for a modified functional:

  \medskip

\begin{figure}[!htb]
 \centering
 \includegraphics[scale=.3]{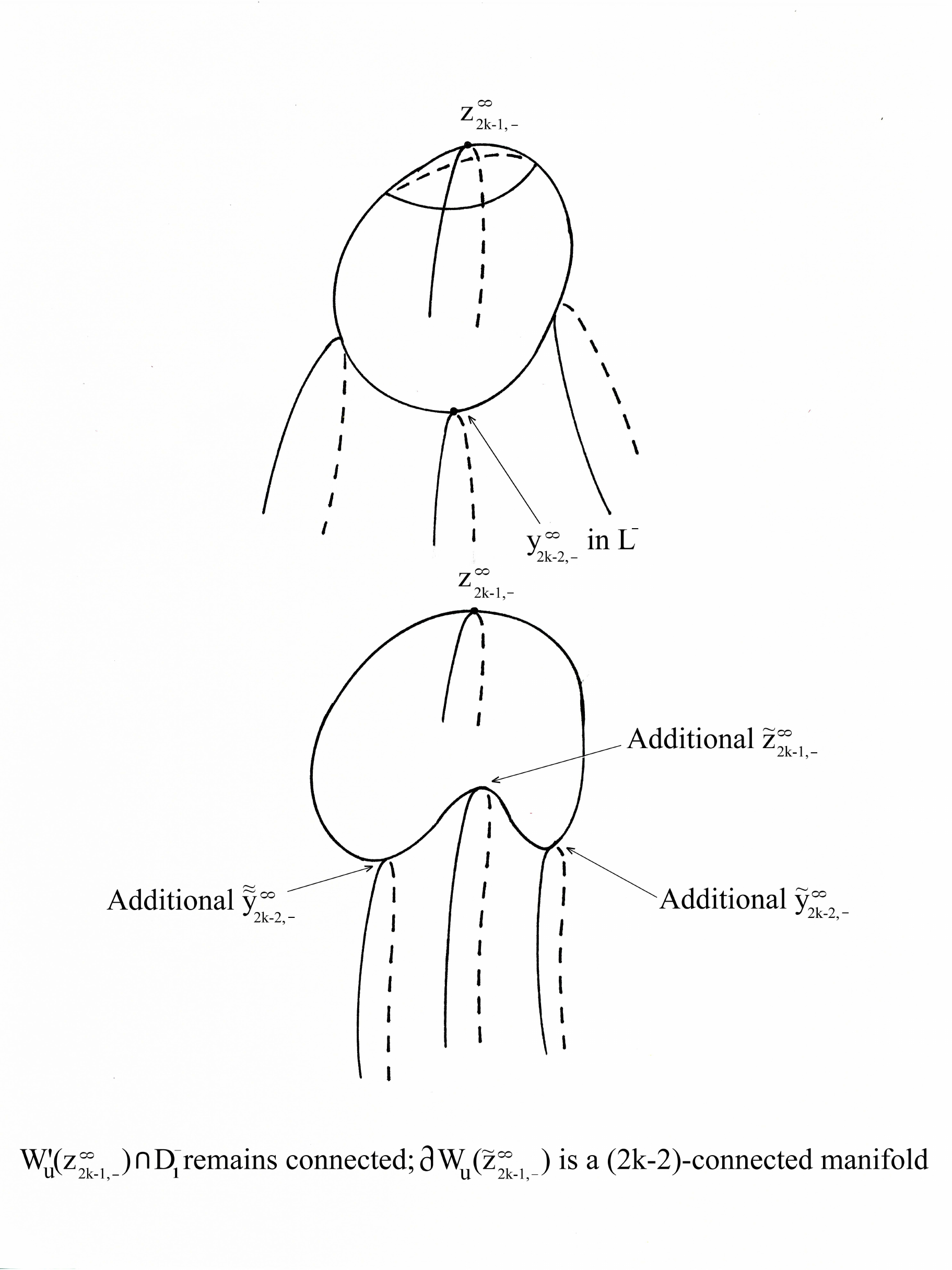}
 \end{figure}


\medskip

  It is important to note that the bottom set for the
  modified $W'_u(z_{2k-1, -}^\infty)$, $W'_u(z_{2k-1, -}^\infty) \cap
  D_1^-$ remains connected since there are only points in the
  unstable sphere of $z_{2k-1, -}^\infty$ which are attracted to the
  critical points at infinity of $L^-$ of index $(2k-2)$. The
  contradiction argument above can therefore run, unchanged.

  \medskip

\noindent{\it  Deleting neighborhoods of periodic orbits in
$c_{2k-1}$}
\medskip

For each periodic orbit $z_i$ dominated by $c_{2k-1}$, w choose a
neighborhood $W_i$ which we delete from $c_{2k-1}$. Using
Proposition 7.24, p608 of [6], which provides an understanding for
the behavior of the flow-lines of $c_{2k-1}$ near $z_i$, we see that
the "$b$"-map of [5] is valued on $\partial W_i \cap c_{2k-1}$ in
$\mathbb{PC}^{k-1} \times \{-1,1\} \cup \mathbb{PC}^{k-2} \times
[-1,1]$. We therefore delete in the pairs of section 2 the $W_i$s
from the first sets of our pairs and we add the $\partial W_i$s to
the second sets of the pairs, leaving the reasoning and the
arguments unchanged.

\bigskip

  \medskip
\noindent{\bf 11. The proof of Theorem 1.3 (i) of [1], of Theorem 1
of the present paper and the proof of the Weinstein Conjecture on
$S^3$, "in the large"}
\medskip

We recall that we have modified in section our functional into the
functional $\tilde J$. $\tilde J^{-1}(\epsilon)$ is $D_1^+ \cup
D_1^-$.  $L^+$ and $L^-$ are to be thought in what follows as small
attracting (for the decreasing pseudo-gradient) neighborhoods of
these sets.

From our results in [5], Propositions 4 and 5, we know that the map
"b" of pairs in homology of dimension $(2k-1)$:

$$H_{2k-1}(\overline {W_u(c_{2k-1})\smallsetminus
   (L^+\cup L^-)},$$ $$(\overline {W_u(c_{2k-1})\smallsetminus
   (L^+\cup L^-)}) \cap [(\partial L^+
    \cup \partial L^-)\cup \tilde {J}_\infty^{-1}(\epsilon)]\cup \overline{\partial_\infty(c_{2k-1}\smallsetminus (L^+ \cup L^-))})\overset{\mathrm{"b"_*}} \longrightarrow$$
   $$ H_{2k-1}(\mathbb{P}\mathbb{C}^{k-1}
   \times [-1,1],\mathbb{P}\mathbb{C}^{k-2}
   \times [-1,1]\cup \mathbb{P}\mathbb{C}^{k-1}
   \times \{-1,[0,1]\})$$

   is onto.

   On the other hand, we know that we have the excision isomorphism (also between pairs):

$$ H_{2k-1}(\overline {W_u(c_{2k-1})\smallsetminus
   (L^+\cup L^-)},$$ $$(\overline {W_u(c_{2k-1})\smallsetminus
   (L^+\cup L^-)})
    \cap [(\partial L^+ \cup \partial L^-)\cup \tilde {J}_\infty^{-1}(\epsilon)]\cup
    (\overline {\partial_\infty(c_{2k-1}\smallsetminus (L^+ \cup L^-)))})\overset{\mathrm{exc}}{\cong}$$



  $$H_{2k-1}(\overline {W_u(c_{2k-1}+h_{2k-1, \infty})\smallsetminus
   (L^\pm)},$$ $$(\overline {W_u(c_{2k-1}+h_{2k-1, \infty})\smallsetminus (L^\pm)}) \cap
    [(\partial (L^\pm)\cup \tilde {J}_\infty^{-1}(\epsilon))]\cup (\overline {W_u(h_{2k-1, \infty}\smallsetminus (L^\pm)})))$$


   We consider the map "$b$", appropriately modified as indicated above.
    We know-this is a key point-that this map extends as an equivariant map to
     $\overline {W_u(h_{2k-1, \infty})\smallsetminus (L^+ \cup L^-)}$ and
     that the restriction of the extension to this set is valued into $\mathbb{PC}^{k-2}\times [-1,1]$.
     We modify slightly our pairs above with the introduction, in the second sets of the pairs, of the additional set
     $B_0$ of section. $J$ is modified into $\tilde J$, the set $\tilde J^{-1}(\epsilon) \cup B_0$
     is alternatively $D_1^+ \cup D_1^- \cup W_u(x^{1, \infty}_-)$. We then find the two pairs of sets $(A,B)$ and $(C,D)$, where:
   $$A=\overline {W_u(c_{2k-1}+h_{2k-1, \infty})\smallsetminus
   (L^+\cup L^-)}$$
    $$B=\overline {W_u(c_{2k-1}+h_{2k-1, \infty})\smallsetminus (L^+\cup
   L^-)}) \cap [(\partial L^+ \cup \partial L^-)\cup \tilde {J}_\infty^{-1}(\epsilon)\cup B_0]$$

   $$C=\overline {W_u(c_{2k-1})\smallsetminus
   (L^+\cup L^-)}$$
   $$D= (\overline {W_u(c_{2k-1})\smallsetminus
   (L^+\cup L^-)}) \cap [(\partial L^+ \cup \partial L^-)\cup \tilde {J}_\infty^{-1}(\epsilon)\cup B_0]\cup
   (\overline {\partial_\infty(c_{2k-1}\smallsetminus (L^+ \cup L^-)))}$$
    The homomorphism:

   $$ H_{2k-1}(A,B) \overset {n_*} \to{\longrightarrow} H_{2k-1}(C,D)$$
   is onto. This follows from the fact that the excision homomorphism above is onto and from the fact that $c_{2k-1}$ is
   assumed to be a minimal cycle (see [1]) of $\partial_{per}$, ie we assume that $c_{2k-1}$ cannot be decomposed
   into smaller cycles for $\partial_{per}$.
   Observe that $A$ is a cycle of dimension $(2k-1)$ relative to $B$, this follows from the relation $(*)$ which we assume to hold.

   Let us also consider the three following maps:

   $$ H_{2k-1}(A,B)\overset {l_*} {\mathrm {\longrightarrow}}  H_{2k-1}(\mathbb{P}\mathbb{C}^{k-1}
   \times [-1,1],\mathbb{P}\mathbb{C}^{r}
   \times [-1,1]\cup \mathbb{P}\mathbb{C}^{k-1}
   \times \{-1,[0,1]\})$$

    $$H_{2k-1}(C,D)\overset {"b"_*} {\mathrm {\longrightarrow}} H_{2k-1}(\mathbb{P}\mathbb{C}^{k-1}
   \times [-1,1],\mathbb{P}\mathbb{C}^{k-2}
   \times [-1,1]\cup \mathbb{P}\mathbb{C}^{k-1}
   \times \{-1,[0,1]\})$$

   $$ H_{2k-1}(\mathbb{P}\mathbb{C}^{k-1}
   \times [-1,1], \mathbb{P}\mathbb{C}^{k-1}
   \times \{-1,[0,1]\}\cup \mathbb{PC}^{r}\times [-1,1])\overset {m_*} {\mathrm {\longrightarrow}}$$

   $$  H_{2k-1}(\mathbb{P}\mathbb{C}^{k-1}
   \times [-1,1],\mathbb{P}\mathbb{C}^{k-2}
   \times [-1,1]\cup \mathbb{P}\mathbb{C}^{k-1}
   \times \{-1,[0,1]\})$$

   The two homomorphisms above $m_*$ and "$b_*$" are onto in dimension $(2k-1)$ (the addition of $D_1^+$
   in the second factor of the pairs $(A,B)$ and $(C,D)$ does not change much to the surjectivity of $"b_*"$ since $U_1$
   maps into a fixed $\mathbb{PC}^r \times [-1,1]$) and the commutation relation $ "b"_* \circ n_*=m_* \circ l_*$ holds.
   It follows that $l_*$ is non-zero. On the other hand, we have the inclusion map

   $$i: (A,B) \overset {i} {\mathrm {\longrightarrow}} (C_\beta \smallsetminus (L^+\cup L^-), (C_\beta -(L^+ \cup L^-)) \cap
   (\partial (L^+ \cup L^-) \cup \tilde {J}^{-1}(\epsilon)\cup B_0)$$

   The map "$b$" extends then in a natural way (this requires the use of general position in order to remove the periodic orbits,
   also the equivariance of the map is as above, on compact sets, with a $p$ in the $e^{ip\tau}$
   that may tend to $\infty$ with the compact sets getting larger, also appropriate powers are taken) into a map:
   $$(C_\beta \smallsetminus (L^+\cup L^-), (C_\beta -(L^+ \cup L^-)) \cap (\partial (L^+ \cup L^-) \cup \tilde {J}^{-1}(\epsilon)\cup B_0))$$
   $$\longrightarrow (\mathbb{P}\mathbb{C}^{\infty}
   \times [-1,1], \mathbb{P}\mathbb{C}^{\infty}
   \times \{-1,[0,1]\}\cup \mathbb{PC}^{r} \times [-1,1])$$

   This implies that $(\overline {W_u(c_{2k-1}+h_{2k-1, \infty})\smallsetminus
   (L^+\cup L^-)},
    (\overline {W_u(c_{2k-1}+h_{2k-1, \infty})\smallsetminus (L^+\cup
   L^-)}) \cap [(\partial L^+ \cup \partial L^-)\cup \tilde {J}_\infty^{-1}(\epsilon)\cup B_0])$
   is not a boundary in $(C_\beta \smallsetminus (L^+\cup L^-), (C_\beta -(L^+ \cup L^-)) \cap [\partial (L^+ \cup L^-) \cup \tilde {J}^{-1}(\epsilon)\cup B_0])$,
   that is that the relation:

   $$\partial c_{2k}^{(\infty)}=c_{2k-1}+h_{2k-1,\infty}$$
   is not possible. The argument is complete.\\

 \noindent{\bf 12. Existence Argument without the basic
assumption}
\medskip

Along a deformation of contact forms, $L^+$ and $L^-$ might change
with the addition or substraction of critical points at infinity
$z^\infty_j$ of index $j$, typically of index $(2k-1)$. The Morse
complex of eg $L^+$ then changes with the addition or the
substraction of a smaller Morse complex. Using the arguments of
Lemma 3, section 6, this smaller Morse complex maps through the
"global" equivariant map "b", see section above, into
$\mathbb{PC}^\infty \times [0,1] \cup \mathbb{PC}^r \times [-1,1]$,
$r$ small when compared to $j$ or $k$. The target value of the
classifying map $l_*$ of section 11 is then unchanged.

The conclusion is that, either using these equivariant/linking
classes, we find a periodic orbit (maybe an iterate) of index
$(2k-1)$, for $k$ large; or there is a periodic orbit of index $1$
connecting $ L^+$ and $L^-$. If there is no such periodic orbit and
these latter sets are connected directly by a critical point at
infinity of index $1$, then, after some reasoning, we find that we
can complete tangencies with other critical points of index $1$
connecting $J_0^{-1}(\epsilon)$ and each of these two sets (we might
need to re-parametrize the flow-lines as in J.Milnor [11], Theorem
4.1 ,pp 37-38, thereby modifying the functional but not the
flow-lines) and completely disconnect these two sets. The existence
argument then proceeds "a la P.Rabinowitz [12]".

To a certain extent, the arguments of this paper indicate that
either we can use the existence argument of H.Hofer [10] and find a
periodic orbit of index $1$ or the equivariant/linking argument of
P.Rabinowitz [12] can be used, one line of proof excluding the other
one. Of course, this is only an indication and not a proof of a
rigorous statement.

 \end{document}